\documentclass[reqno]{amsart}
\usepackage{latexsym,amsfonts,amssymb,amsmath,euscript,amsthm}
\usepackage{graphicx} 
\usepackage{hyperref}         
\usepackage{color,fancybox}

\newcommand{\cqd}{\hspace{10pt}\fbox{}}
\newcommand{\R} {\mathbb{R}}
\newcommand{\N} {\mathbb{N}}

\newcommand{\eps}{\epsilon}

\def\eto{\buildrel \epsilon\to 0\over\longrightarrow }

\begin{document}
\title[Homogenization in thin domains]{Homogenization in a thin domain with an \\
oscillatory boundary}
\author[J. M. Arrieta]{Jos\'{e} M. Arrieta$^{*, \dag}$}
\thanks{$^*$ Corresponding author:  Jos\'e M. Arrieta,  Departamento de Matem\'atica Aplicada, Facultad de Matem\'aticas, 
Universidad Complutense de Madrid, 28040  Madrid, Spain. e-mail: arrieta@mat.ucm.es}
\thanks{$^\dag$ Partially
supported by:  MTM2009-07540  DGES, Spain; PHB2006-003 PC and PR2009-0027 from MICINN;
     and  GR58/08, Grupo 920894 (BSCH-UCM, Spain)}
\address[J. M. Arrieta]{Departamento de Matem\'atica Aplicada,
Facultad de Ma\-te\-m\'a\-ti\-cas, Universidad Complutense de
Madrid, 28040 Madrid, Spain.}  \email{arrieta@mat.ucm.es}

\author[M.C.Pereira]{Marcone C. Pereira$^\ddag$}\thanks{$^\ddag$Partially
supported by FAPESP 2008/53094-4, CAPES DGU 127/07 and CNPq
305210/2008-4.}
\address[M. C. Pereira]{Escola de Artes, Ci\^encias e Humanidades \\ Universidade de S\~ao
Paulo, S\~ao
Paulo, SP, Brazil}   \email{marcone@usp.br}

\date{}

\begin{abstract}
In this paper we analyze the behavior of the Laplace operator with Neumann boundary 
conditions in a thin domain of the type $R^\epsilon = \{ (x_1,x_2) \in \R^2 \; | \;  x_1 \in (0,1), \, 0 < x_2 < \epsilon \, G(x_1, x_1/\eps) \}
$ where the function $G(x,y)$ is periodic in $y$ of period $L$.  Observe that the upper boundary of the thin domain presents a highly oscillatory behavior and, moreover, the height of the thin domain, the amplitude
and period of the oscillations are all of the same order, given by the small parameter $\epsilon$.

\par\bigskip \noindent R\'esum\'e: Dans cet article, nous analysons le comportement de l'op\'erateur de Laplace avec conditions aux limites de Neumann 
dans un domaine fine du type $R^\epsilon = \{ (x_1,x_2) \in \R^2 \; | \;  x_1 \in (0,1), \, 0 < x_2 < \epsilon \, G(x_1, x_1/\eps) \}
$
lorsque la fonction $G(x,y)$ est p\'eriodique dans la variable $y$ de p\'eriode $L$. 
On observe que la limite sup\'erieure du domaine fine pr\'esente une comportement hautement oscillatoire et, en outre, l' hauteur du domaine, l'amplitude et la p\'eriode des oscillations sont tous du m\^eme ordre, donn\'e par un petit param\`etre  $\epsilon$.


\end{abstract}

\maketitle


\numberwithin{equation}{section}
\newtheorem{theorem}{Theorem}[section]
\newtheorem{lemma}[theorem]{Lemma}
\newtheorem{corollary}[theorem]{Corollary}
\newtheorem{proposition}[theorem]{Proposition}
\newtheorem{remark}[theorem]{Remark}
\allowdisplaybreaks

\section{Introduction}

In this work we study the asymptotic behavior of the solutions of the Neumann problem for the Laplace operator
\begin{equation} \label{OP0}
\left\{
\begin{gathered}
- \Delta w^\epsilon + w^\epsilon = f^\epsilon
\quad \textrm{ in } R^\epsilon \\
\frac{\partial w^\epsilon}{\partial \nu ^\epsilon} = 0
\quad \textrm{ on } \partial R^\epsilon
\end{gathered}
\right. 
\end{equation}
with $f^\epsilon \in L^2(R^\epsilon)$,  $\nu^\epsilon = (\nu^\epsilon_1, \nu^\epsilon_2)$ is the unit outward normal to $\partial R^\epsilon$
and $\frac{\partial }{\partial \nu^\epsilon}$ is the outside normal derivative.  The domain $R^\epsilon$ is a thin domain with a highly oscillating boundary which is given by
$$
R^\epsilon = \{ (x_1,x_2) \in \R^2 \; | \;  x_1 \in (0,1),  \quad
 0 < x_2 < \epsilon \, G_\epsilon(x_1) \}
$$
where $G_\eps(\cdot)$ is a function satisfying $0<G_0\leq G_\eps(\cdot)\leq G_1$ for some 
fixed positive constants $G_0$, $G_1$ and such that oscillates as the parameter $\eps\to 0$. We may think, for instance,  that the function $G_\eps$ is of the form $G_\eps(x)=a(x)+b(x)g(x/\eps^\alpha)$ where $a, b:I \mapsto \R$ are piecewise $\mathcal{C}^1$-functions defined on $I = (0,1)$,  $g:\R\to \R$ is an $L$-periodic smooth function and $\alpha\geq 0$, see Figure 1. This includes the case where the function $G_\eps(\cdot)$ is a purely periodic function, for instance, $G_\eps(x)=2+\sin(x/\eps^\alpha)$ but also includes the 
 case where the function $G_\eps$ is not periodic and the amplitude is modulated by a function. As a matter of
 fact, we will be able to treat more general cases, see hypothesis (H) below, but to stay the general ideas in the
 introduction we may consider the proptotype function  $G_\eps(x)=a(x)+b(x)g(x/\eps^\alpha)$.

\begin{figure}[h]
\centering{\includegraphics[width=9cm,height=3.2cm]{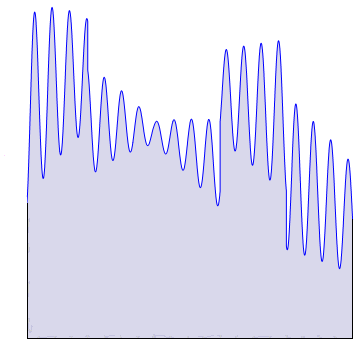}}
\label{thin-domain}
\caption{The thin domain $R^\eps$}
\end{figure}

The existence and uniqueness of solutions for problem (\ref{OP0}) for each $\eps>0$ is garanteed by Lax-Milgram Theorem.
We are interested here in analyzing the behavior of the solutions as $\eps\to 0$, that is, as the domain gets
thinner and thinner although with a high oscillating boundary. 

Observe that the domain is thin since $R^\epsilon\subset (0,1)\times (0, \eps G_1)$ and its upper boundary oscillates due to the function $g(x/\eps^\alpha)$, (if $\alpha>0$ and $g$ is not a constant function).

Since the domain $R^\eps$ is thin, ``approaching'' the line segment $(0,1)\subset \R$, it is reasonable
to expect that the family of solutions  $w^\eps$ will converge to a function of just one variable and that
this function will satisfy an equation of the same type of (\ref{OP0}) but in one dimension, say $Lu+u=h$ in $(0,1)$ with
some boundary conditions,  where
$L$ is a second order elliptic operator in one dimension.  As a matter of fact, if the function $G_\eps(x)$ is independent
of $\eps$, (say $\alpha=0$ or $g\equiv 0$),  that is, the thin domain does not present any oscillations whatsoever:
$$R^\eps=\{ (x_1,x_2):   0<x_1<1, 0<x_2<\eps G(x)\}$$
the limit equation is given by 

\begin{equation}\label{thin-domain-hale-raugel}
\left\{
\begin{gathered}
- \frac{1}{G(x)}(G(x)w_x)_x + w(x) = f(x)\quad  x\in (0,1)\\
w_x(0)=w_x(1)=0
\end{gathered}
\right. 
\end{equation}
see for instance \cite{HR, Raugel}.  Observe that the geometry of the thin domain enters into the limit equation
through the diffusion coefficient.

Moreover, if we consider $0\leq \alpha<1$ 
and if we assume that $a(x)+b(x)g(x/\epsilon^\alpha)\to h(x)$ w-$L^2(0,1)$ and $\frac{1}{a(x)+b(x)g(x/\epsilon^\alpha)}\to k(x)$ w-$L^2(0,1)$
(observe that $h(x)k(x)\geq 1$ a.e. and in general it is not true that $h(x)k(x)\equiv 1$), then the limit problem is 
$$
\left\{
\begin{gathered}
-\frac{1}{h(x)} \left(\frac{1}{k(x)} v_x \right)_x + v = f, \hbox{ in }(0,1)\\
v_x(0)=v_x(1)=0
\end{gathered}
\right. 
$$
see \cite{JA}.    Observe that this case contains the previous one. If $\alpha=0$, then $h(x)=a(x)+b(x)g(x)\equiv G(x)$ and $k(x)=\frac{1}{a(x)+b(x)g(x)}=\frac{1}{G(x)}$, and we recover equation (\ref{thin-domain-hale-raugel}).

\par\bigskip
In this work we are interested in addressing the case $\alpha=1$, that is $G_\eps(x)=a(x)+b(x)g(x/\epsilon)$, 
where none of the techniques used to solve the previous ones apply.  Observe that this situation  is very resonant: the height of the domain, the amplitude of the oscillations at the boundary and the period of the oscillations are of the same order $\eps$.  Moreover we are interested in addressing not only the purely periodic case, that is, the case where the
function $G_\eps(x)=G(x/\eps)$ for some $L$-periodic smooth function $G$ but also the general case where the
amplitude of the oscillation depend on $x$ in a continuous fashion, that is, in our prototype case, the situation 
where $a$ and $b$ are not piecewise constant, but piecewise continuous function.

The purely periodic case can be addressed by somehow standard techniques in homogenization theory, as developed 
in \cite{BLP, CP, SP}. If $G_\eps(x)=G(x/\eps)$ where $G$ is an $L$-periodic $C^1$ function and if we denote by
$$
Y^* = \{ (y_1,y_2) \in \R^2 \; : \; 0< y_1 < L, \quad 0 < y_2 < G(y_1) \}
$$
then the limit equation is shown to be 
\begin{equation} \label{GLP}
\left\{
\begin{gathered}
-q_0w_{xx} + w =  f(x), \quad x \in (0,1)\\
w'(0)=w'(1)=0
\end{gathered}
\right.
\end{equation} 
where 
$$
\begin{gathered}
q_0=\frac{1}{|Y^*|} \int_{Y^*} \Big\{ 1 - \frac{\partial X}{\partial y_1}(y_1,y_2) \Big\} dy_1 dy_2, \qquad 
\end{gathered}
$$ 
and $X$ is the unique solution of 
\begin{equation} \label{AUXI}
\left\{
\begin{array}{l}
- \Delta X  =  0  \textrm{ in } Y^*  \\
\displaystyle\frac{\partial X}{\partial N}  =  0  \textrm{ on } B_2  \\
\displaystyle\frac{\partial X}{\partial N}  =  N_1 \textrm{ on } B_1  \\
X(0,y_2)  =  X(L,y_2) \textrm{ on } B_0  \\
\displaystyle\int_{Y^*} X \; dy_1 dy_2  =  0 
\end{array}
\right.
\end{equation}
where $B_0$ is the lateral boundary, $B_1$ is the upper boundary and $B_2$ is the lower
boundary of $\partial Y^*$, that is
$$B_0=\{(0,y_2): 0<y_2<G(0)\}\cup \{(L,y_2): 0<y_2<G(L)\}$$ 
$$B_1=\{ (y_1, G(y_1)):  0<y_1<L\}$$
$$B_2=\{ (y_1, 0):  0<y_1<L\}.$$
We refer to \cite{ACPS} for a complete analysis of the purely periodic case of a nonlinear parabolic problem. 

If we want to analyze now the case where the function $G_\eps$ is given by $G_\eps(x)=a(x)+b(x)g(x/\eps)$ and the functions $a$, $b$
are smooth but not necessarily constant, it is intuitively true that the limiting equation should behave like (\ref{GLP}) with a diffusion coefficient $q_0$ depending on $x$ somehow.  Actually if we look at the thin domain in a small neighborhood of a point $\xi\in (0,1)$, we will approximately see a domain with very high oscillations but locally the domain behaves like the thin domain with the function $x\to a(\xi)+b(\xi)g(x/\eps)$. Therefore, it is expected that if we freeze the coefficients of the limit equation in a fixed point $\xi\in (0,1)$ 
we should recover equation (\ref{GLP}).   Although this intuitive argument will turn out to be true, it does not give us a complete information about the limit equation, specially, the way in which the dependence on $x$ is explicitly  stated in the limit equation. For instance, it is not clear at this stage whether the limit equation should be
$-(q(x)w_x)_x+w=f$ or $-q(x) w_{xx}+w=f$ with $q(x)=\frac{1}{|Y^*(x)|} \int_{Y^*(x)} \{ 1 - \frac{\partial X(x)}{\partial y_1}(y_1,y_2) \} dy_1 dy_2,$ or maybe $-\frac{1}{|Y^*(x)|}(r(x)w_x)_x+w=f$ where $r(x)= \int_{Y^*(x)}\{ 1 - \frac{\partial X(x)}{\partial y_1}(y_1,y_2) \} dy_1 dy_2,$ or maybe other. Observe that all these equations coincide if we consider the
purely periodic case.   

\par\medskip 

In order to accomplish our goal and obtain the limit equation in the general case, we propose a technique that consists in solving first the piecewise periodic case, that is, the case where the functions $a(x)$ and $b(x)$ are piecewise constant and then 
do an approximation argument to obtain the limit equation in the general case.  This is a subtle argument since we are
approximating first the functions $a$ and $b$ by piecewise constant functions, say $a^\delta(x)$ , $b^\delta(x)$ 
so that  $|a(x)-a^\delta(x)|+|b(x)-b^\delta(x)|\leq \delta$ and obtain the limit equation for $\delta>0$ fixed, passing to the limit as $\eps\to 0$. Then, in this limit equation, which will depend on $\delta$, we pass to the limit as $\delta\to 0$. 
 But the limit we are interested in is taking first $\delta\to 0$ for $\eps>0$ fixed, so we obtain the domain given by
the function $a(x)+b(x)g(x/\eps)$, and then we pass to the limit as $\eps\to 0$.  But,   a priori there is
no garantee that these two limits commute.  We will actually show that these two limits commute by proving that
the solutions of problem (\ref{OP0}) in two domains $\Omega^\eps=\{(x,y): 0<x<1, 0<y<G_\eps(x)=a(x)+b(x)g(x/\eps)\}$
and  $\tilde \Omega^\eps=\{(x,y): 0<x<1, 0<y<\tilde G_\eps(x)=\tilde a(x)+\tilde b(x)g(x/\eps)\}$, are close in the $H^1$ norm uniformly in $\eps$ when $a,b$ and $\tilde a, \tilde b$ are close. This result, which can be regarded as a domain perturbation result but uniformly in $\eps$,  guarantee that the two limits commute and 
will show that the limit problem is given as above.  We remark that this domain perturbation result is not
deduced from standard and known results on domain perturbation techniques since we are able to compare the solutions of an elliptic problem in two families of domains which also depend on a parameter and the way
this domains depend on $\eps$ is not smooth at all. 

We strongly believe that this technique of solving first the piecewise periodic case and then use an approximation
argument, {\sl uniform in the parameter $\epsilon$}, can be used in other problems to obtain the appropriate 
homogenized limit for the non periodic case. 

\par\medskip 

Following this agenda, we solve first the piecewise periodic case, that is, the case where the
function $G_\eps(x)=a^\delta(x)+b^\delta(x)g(x/\eps)$ and $a^\delta$, $b^\delta$ are piecewise constant. We consider
the points $0=\xi_0<\xi_1<\ldots<\xi_{N-1}<\xi_N=1$ and assume that the functions $a^\delta$ and $b^\delta$ are constant in each
of the interval $(\xi_{i-1},\xi_i)$, say $a^\delta(x)=a_i$, $b(x)=b_i$. We show that the limit equation is of the same form (\ref{GLP})
in each of the intervals $(\xi_{i-1},\xi_i)$, that is, 
\begin{equation} \label{GLP-piecewise}
-q_iw_{xx} + w =  f(x), \quad x \in (\xi_{i-1},\xi_i), \quad, i=1,\ldots, N
\end{equation} 
where 
$$
\begin{gathered}
q_i=\frac{1}{|Y^*_i|} \int_{Y^*_i} \Big\{ 1 - \frac{\partial X_i}{\partial y_1}(y_1,y_2) \Big\} dy_1 dy_2, \qquad 
\end{gathered}
$$ 
and $X_i$ is the unique solution of (\ref{AUXI}) in the cell $Y^*_i$ where 
$$
Y^*_i= \{ (y_1,y_2) \in \R^2 \; : \; 0< y_1 < L, \quad 0 < y_2 < a_i+b_i g(y_1) \}.
$$
Moreover, equation (\ref{GLP-piecewise}) is suplemented with Neumann boundary conditions at
$x=0$, $x=1$ and with some ``matching'' condition at the points $\xi_{i}$, $i=1,\ldots, N-1$, which
are continuity of the function and some kind of Kirchoff type conditions, see Section \ref{piecewise-periodic} for details. 
Actually, if we look at the variational formulation of the limit equation, we obtain 
\begin{equation} \label{VFPDL-piecewise}
\int_{I} \Big\{ q^\delta(x) \, p^\delta(x) \; w_x(x) \, \varphi_x(x) 
+ p^\delta(x) \, w(x) \, \varphi(x) \Big\} dx = \int_{I} p^\delta(x) f(x) \, \varphi \, dx
\end{equation}
where $q^\delta(x)=q_i$ and $p^\delta(x)=\frac{|Y^*_i|}{L}$  in $(\xi_{i-1},\xi_i)$.

Now we will be able to pass to the limit in this equation as $\delta \to 0$ and obtain the limit problem:

\begin{equation} \label{GLP-cont}
\left\{
\begin{gathered}
-\frac{1}{p(x)}(r(x)w_{x})_x + w =  f(x), \quad x \in (0,1)\\
w'(0)=w'(1)=0
\end{gathered}
\right.
\end{equation} 
where 
$$
\begin{gathered}
r(x)= p(x)q(x)=\frac{1}{L}\int_{Y^*(x)} \Big\{ 1 - \frac{\partial X(x)}{\partial y_1}(y_1,y_2) \Big\} dy_1 dy_2, \qquad 
\end{gathered}
$$ 
$$p(x)=\frac{|Y^*(x)|}{L}$$
and $X(x)$ is the unique solution of (\ref{AUXI}) in the basic cell $Y^*(x)$ which depends on the variable $x$ and 
it is given by 
$$
Y^*(x)= \{ (y_1,y_2) \in \R^2 \; : \; 0< y_1 < L, \quad 0 < y_2 < a(x)+b(x)g(y_1) \}.
$$
Equation (\ref{GLP-cont}) is the limit equation we were looking for.  

\par\bigskip 

Finally, we would like to remark that although we will treat the Neumann boundary condition problem,  
we may also impose different conditions in the lateral boundaries of the thin
domain $R^\eps$, while preserving the Neumann type boundary condition in the upper and lower boundary.  
That is for problem (\ref{OP0}) we may consider conditions of the Dirichlet
type, $w^\eps=0$, or even Robin, $\frac{\partial w^\eps}{\partial N}+\beta w^\eps=0$ in the lateral
boundaries $\{ (0,y): 0<y<\eps G_\eps(0)\}$,  $\{ (1,y): 0<y<\eps G_\eps(1)\}$.
The limit problem will preserve this boundary condition. That is, in problem (\ref{GLP-cont}) we will
obtain the conditions $w=0$ or $\frac{\partial w}{\partial N}+\beta w=0$ at $x=0$ and $x=1$.

\par\bigskip 

We describe the contents of the paper. In Section \ref{PRE} we set up the notation, obtain some technical 
results that will be needed in the proofs and state the main result.   In Section \ref{piecewise-periodic} we 
obtain the result for the piecewise periodic case. In Section \ref{BS} we show the continuous dependence result
on the domain, uniform in $\epsilon$, that will be the key argument to obtain the limit problem.  In Section \ref{general-case} we provide a proof of the main result. Finally, we include an Appendix where we analyze the behavior
of the basic function $X$ solution of (\ref{AUXI}) as we perturb $G$.

\bigskip

\par\noindent {\bf Acknowledgments.} Part of this work was done while the second author was visiting the Applied Math Department at Complutense University in Madrid. He express his kind gratitude to the Department.  

\noindent We would also like to thank E. Zuazua, C. Castro, M. Eugenia P\'erez, M. Lobo and D. G\'omez for the all the discussion and good suggestions on  this problem.  

\section{Basic facts, notation and main result} \label{PRE}

We consider the one parameter family of functions $G_\eps: I \to (0,\infty)$, where
$I=(0,1)$, $\eps\in (0, \eps_0)$ for some $\eps_0>0$. 
We will assume the following hipothesis

\begin{itemize}
\item[(H)] There exist two positive constants $G_0$, $G_1$ such that 
\begin{equation} \label{HG}
\begin{gathered}
0< G_0 \le G_\epsilon(x) \le G_1, \quad  \forall x\in I, \quad \forall \eps\in (0,\eps_0).
\end{gathered}
\end{equation}
Moreover, the  functions $G_\eps(\cdot)$ are of 
the type $G_\eps(x)=G(x,x/\eps)$, where  the function
\begin{equation}\label{def-G}
\begin{array}{rl}
G:I\times \R &\longrightarrow (0,+\infty) \\
 (x,y)&\longrightarrow G(x,y)
 \end{array}
 \end{equation}
 is $L$-periodic in the second variable, that is,  $G(x,y+L)=G(x,y)$ and piecewise $C^1$ with
 respect to the first variable, that is, there exists a finite number of points $0=\xi_0<\xi_1<\cdots<\xi_{N-1}<\xi_N=1$ such
 that the function $G:(\xi_i,\xi_{i+1})\times \R\to (0,+\infty)$ is $C^1$ and  such that $G$, $G_x$ and  $G_y$ are uniformly bounded in $(\xi_i,\xi_{i+1})\times \R$ and have limits when we approach $\xi_i$ and $\xi_{i+1}$. 
\end{itemize}

One important example of a function satisfying the above conditions is  
$$G(x,y)=a(x)+\sum_{r=1}^N b_r(x)g_r(y)$$
where the functions $a$, $b_1$,..,$b_N$ are piecewise $C^1$ in $I=(0,1)$ and the functions
$g_1$,..,$g_N$ are all $C^1$ and $L$-periodic.

We define the thin domain as
$$
R^\epsilon = \{ (x_1,x_2) \in \R^2 \; | \;  x_1 \in I,  \quad
 0 < x_2 < \epsilon \, G_\epsilon(x_1) \}.
$$

In this work we study the asymptotic behavior of the solutions of the Neumann problem for the Laplace operator
\begin{equation} \label{OP}
\left\{
\begin{gathered}
- \Delta w^\epsilon + w^\epsilon = f^\epsilon
\quad \textrm{ in } R^\epsilon \\
\frac{\partial w^\epsilon}{\partial \nu ^\epsilon} = 0
\quad \textrm{ on } \partial R^\epsilon
\end{gathered}
\right. 
\end{equation}
with $f^\epsilon \in L^2(R^\epsilon)$ and where $\nu^\epsilon = (\nu^\epsilon_1, \nu^\epsilon_2)$ is the unit outward normal to $\partial R^\epsilon$
and $\frac{\partial }{\partial \nu^\epsilon}$ is the outside normal derivative. 

From Lax-Milgran Theorem, we have that  problem (\ref{OP}) has a unique solution for each $\epsilon > 0$. 
We are interested here in analyzing the behavior of the solutions as $\eps\to 0$, that is, as the domain gets
thinner and thinner although with a high oscillating boundary.

To study the convergence of (\ref{OP}) in the thin oscillating domain $R^\epsilon$, 
we consider the equivalent linear elliptic problem
\begin{equation} \label{P}
\left\{
\begin{gathered}
 - \frac{\partial^2 u^\epsilon}{{\partial x_1}^2} - 
\frac{1}{\epsilon^2} \frac{\partial^2 u^\epsilon}{{\partial x_2}^2} + u^\epsilon = f^\epsilon
\quad \textrm{ in } \Omega^\epsilon \\
\frac{\partial u^\epsilon}{\partial x_1} N_1^\epsilon + \frac{1}{\epsilon^2} \frac{\partial u^\epsilon}{\partial x_2}N_2^\epsilon = 0
\quad \textrm{ on } \partial \Omega^\epsilon
\end{gathered}
\right.
\end{equation}
where $f^\eps\in L^2(\Omega^\eps)$  satisfies 
\begin{equation} \label{ESTF}
\| f^\epsilon \|_{L^2(R^\epsilon)} \le C
\end{equation}
for some $C > 0$ independent of $\epsilon$ and now $N^\eps=(N_1^\eps,N_2^\eps)$ is the outward
unit normal to $\partial\Omega^\eps$ and 
$\Omega^\epsilon \subset \R^2$ is given by
\begin{equation} \label{domain}
\Omega^\epsilon = \{ (x_1,x_2) \in \R^2 \; | \;  x_1 \in I, \quad 0 < x_2 < G_\epsilon(x_1) \}.
\end{equation}

Observe that the equivalence between the problems (\ref{OP}) and (\ref{P}) is obtained by changing the scale 
of the domain $R^\epsilon$ through the transformation $(x,y) \to (x, \epsilon y)$, (see \cite{HR} for more details).
Moreover, the domain $\Omega^\eps$ is not ``thin'' anymore but presents very wild oscillations  at the top boundary, 
although the presence of a high diffusion coefficient in front of the derivative with respect the second variable balance
the effect of the wild oscillations. 

The domain $\Omega^\eps$ is related to the ones analyzed in the papers \cite{BCh, ABMG,DP}
but the fact that in our case we have a very high diffusion in the $y$-direction makes our analysis and results different from these other papers. 

\par\medskip
We write $H^1_\epsilon(U)$ for the space $H^1(U)$ with the equivalent norm
$$
\| w \|_{H^1_\epsilon(U)}^2 = \| w \|_{L^2(U)}^2 + \Big\| \frac{\partial w}{\partial x_1} \Big\|_{L^2(U)}^2
+ \frac{1}{\epsilon^2} \Big\| \frac{\partial w}{\partial x_2} \Big\|_{L^2(U)}^2
$$
given by the inner product 
$$
(\phi, \varphi)_{H^1_\epsilon(U)} = \int_{U} \Big\{ \phi \, \varphi + \frac{\partial \phi}{\partial x_1} \frac{\partial \varphi}{\partial x_1} 
+ \frac{1}{\epsilon^2} \frac{\partial \phi}{\partial x_2} \frac{\partial \varphi}{\partial x_2}
 \Big\} dx_1 dx_2
$$
where $U$ is an arbitrary open set of $\R^2$, which may depend also on $\eps$.

The variational formulation of (\ref{P}) is find $u^\epsilon \in H^1(\Omega^\epsilon)$ such that 
\begin{equation} \label{VFP}
\int_{\Omega^\epsilon} \Big\{ \frac{\partial u^\epsilon}{\partial x_1} \frac{\partial \varphi}{\partial x_1} 
+ \frac{1}{\epsilon^2} \frac{\partial u^\epsilon}{\partial x_2} \frac{\partial \varphi}{\partial x_2}
+ u^\epsilon \varphi \Big\} dx_1 dx_2 = \int_{\Omega^\epsilon} f^\epsilon \varphi dx_1 dx_2, 
\, \forall \varphi \in H^1(\Omega^\epsilon).
\end{equation}
Remark that the solutions $u^\epsilon$ satisfy an uniform a priori estimate on $\epsilon$.
In fact, taking $\varphi = u^\epsilon$ in  expression (\ref{VFP}),  we obtain 
\begin{equation} \label{priori}
\begin{gathered}
\Big\| \frac{\partial u^\epsilon}{\partial x_1} \Big\|_{L^2(\Omega^\epsilon)}^2
+ \frac{1}{\epsilon^2}\Big\| \frac{\partial u^\epsilon}{\partial x_2} \Big\|_{L^2(\Omega^\epsilon)}^2
+ \| u^\epsilon \|_{L^2(\Omega^\epsilon)}^2
\le \| f^\epsilon \|_{L^2(\Omega^\epsilon)} \| u^\epsilon \|_{L^2(\Omega^\epsilon)}.
\end{gathered}
\end{equation}
Consequently, 
\begin{equation} \label{EST0}
\begin{gathered}
\| u^\epsilon \|_{L^2(\Omega^\epsilon)}, \Big\| \frac{\partial u^\epsilon}{\partial x_1} \Big\|_{L^2(\Omega^\epsilon)}
\textrm{ and } \frac{1}{\epsilon} \Big\| \frac{\partial u^\epsilon}{\partial x_2} \Big\|_{L^2(\Omega^\epsilon)} 
\le \| f^\epsilon \|_{L^2(\Omega^\epsilon)} \quad \forall \epsilon > 0.
\end{gathered}
\end{equation}
Hence, it follows from (\ref{ESTF}) that there exists $C > 0$, independent of $\epsilon > 0$, such that
\begin{equation} \label{EST00}
\begin{gathered}
\| u^\epsilon \|_{L^2(\Omega^\epsilon)}, \Big\| \frac{\partial u^\epsilon}{\partial x_1} \Big\|_{L^2(\Omega^\epsilon)}
\textrm{ and } \frac{1}{\epsilon} \Big\| \frac{\partial u^\epsilon}{\partial x_2} \Big\|_{L^2(\Omega^\epsilon)} 
\le C.
\end{gathered}
\end{equation}

Observe that the solutions $u^\epsilon$ of (\ref{P}) can also be characterized as the minimum of a functional.
That is, $u^\epsilon \in H^1(\Omega^\epsilon)$ is the solution of (\ref{P}) if only if 
\begin{equation} \label{FUNC}
V_\epsilon(u^\epsilon) = \min_{\varphi \in H^1(\Omega^\epsilon)} V_\epsilon(\varphi),
\end{equation}
where
$$
\begin{gathered}
V_\epsilon: H^1(\Omega^\epsilon) \mapsto \R \\
V_\epsilon(\varphi) = 
\frac{1}{2} \int_{\Omega^\epsilon} \Big\{ \frac{\partial \varphi}{\partial x_1}^2 
+ \frac{1}{\epsilon^2} \frac{\partial \varphi}{\partial x_2}^2 
+ \varphi^2 \Big\} dx_1 dx_2 - \int_{\Omega^\epsilon} f^\epsilon \varphi dx_1 dx_2.
\end{gathered}
$$
It follows from (\ref{VFP}) with $\varphi = u^\epsilon$ that  
$$
V_\epsilon(u^\epsilon) = - \frac{1}{2} \int_{\Omega^\epsilon} f^\epsilon u^\epsilon dx_1 dx_2.
$$
Hence, due to (\ref{ESTF}) and (\ref{EST0}) we obtain
\begin{equation} \label{ESTV}
|V_\epsilon(u^\epsilon)| \leq \frac{1}{2} \| f^\epsilon \|_{L^2(\Omega^\epsilon)} \| u^\epsilon \|_{L^2(\Omega^\epsilon)} \leq C^2\end{equation}

\par\bigskip

Important tools for the analysis are appropriate extension operators for functions defined in the sets 
$\Omega^\eps$.  We will obtain such operator and we will be able to construct it even for more general domains
that the ones like $\Omega^\eps$. 

Hence, let us consider the following open sets:
$$
\begin{gathered}
\mathcal{O} = \{ (x_1, x_2) \in \R^2 \; | \; x_1 \in I 
\textrm{ and } 0 < x_2 < G_1 \} \\
\mathcal{O}^\epsilon = \{ (x_1, x_2) \in \R^2 \; | \; x_1 \in I 
\textrm{ and } 0 < x_2 < G_\epsilon (x_1) \}
\end{gathered}
$$
where  $I \subset \R$ is an open interval, $G_\epsilon: I \mapsto \R$ is a piecewise $\mathcal{C}^1$-function satisfying 
$0 < G_0 \le G_\epsilon (x_1) \le G_1$ for all $x \in I$ and $\epsilon > 0$ and there exists $0=\xi_0<\xi_1<\ldots<\xi_{N-1}<\xi_N=1$
such that $G_\eps$ is $\mathcal{C}^1$ in the intervals $(\xi_i,\xi_{i+1})$. Let us define 
\begin{equation}\label{definition-eta}
\eta(\epsilon) = \sup_{x \in I\setminus{\xi_0,\ldots,\xi_N}} \{ | G'_\epsilon(x) | \}
\end{equation}
and assume $\eta(\eps)<+\infty$ for fixed $\eps$, although in general $\eta(\eps)\to +\infty$ as $\eps\to 0$. 

Also, we denote by 
$$
\hat {\mathcal{O}} = \mathcal{O}\setminus    \cup_{i=1}^N\{ (\xi_i,x_2):   G_0<x_2<G_1\}.$$

Notice that $\mathcal{O}^\epsilon \subset \mathcal{O}$.

\begin{lemma} \label{EOT}

With the notation above, there exists an extension  operator 
$$
P_{\epsilon} \in \mathcal{L}(L^p(\mathcal{O}^\epsilon),L^p(\hat{\mathcal{O}})) 
  \cap \mathcal{L}(W^{1,p}(\mathcal{O}^\epsilon),W^{1,p}(\hat{\mathcal{O}}))
   \cap \mathcal{L}(W^{1,p}_{\partial_l}(\mathcal{O}^\epsilon),W^{1,p}_{\partial_l}(\hat{\mathcal{O}}))
$$
(where $W^{1,p}_{\partial_l}$ is the set of functions in $W^{1,p}$ which are zero in the lateral boundary $\partial_l$) 
and a constant $K$ independent of $\epsilon$ and $p$ such that
\begin{equation} 
\begin{gathered} \label{EQOP}
\| P_{\epsilon} \varphi \|_{L^p(\hat{\mathcal{O}})} \le K \, \| \varphi \|_{L^p(\mathcal{O}^\epsilon)}  \\
\Big\| \frac{\partial P_{\epsilon} \varphi}{\partial x_1} \Big\|_{L^p(\hat{\mathcal{O}})} 
\le K \, \Big\{ \Big\| \frac{\partial \varphi}{\partial x_1} \Big\|_{L^p(\mathcal{O}^\epsilon)} 
+  \eta(\epsilon) \, \Big\| \frac{\partial \varphi}{\partial x_2} \Big\|_{L^p(\mathcal{O}^\epsilon)} \Big\} \\
\Big\| \frac{\partial P_{\epsilon} \varphi}{\partial x_2} \Big\|_{L^p(\hat{\mathcal{O}})} 
\le K \, \Big\| \frac{\partial \varphi}{\partial x_2} \Big\|_{L^p(\mathcal{O}^\epsilon)} 
\end{gathered}
\end{equation}
for all $\varphi \in W^{1,p}(\mathcal{O}^\epsilon)$ where $1 \le p \le \infty$ and  $\eta(\eps)$ is defined in \eqref{definition-eta}. 
\end{lemma}
\begin{proof}
Observe first that the set $\mathcal{O}^{0}=(0,1)\times (0,G_0)\subset \mathcal{O}^{ \epsilon}$ for 
all $\eps$.   Hence, if we have that $G_1\leq 2G_0$, which implies that $G_\eps(x_1)\leq 2G_0$, we
can define the operator:
$$
(P_\epsilon \varphi) (x_1, x_2) = 
\left\{ 
\begin{array}{ll}
\varphi(x_1, x_2) & (x_1, x_2 ) \in \mathcal{O}^\epsilon \\
\varphi(x_1, 2 G_\epsilon (x_1) - x_2 ) & (x_1, x_2 ) \in \mathcal{O} \slash \mathcal{O}^\epsilon.
\end{array}
\right.
$$

Observe that this operator is obtained through a ``reflection'' procedure in the $x_2$ direction 
along the oscillating boundary.   It is straigth forward to check that this operator satisfies
(\ref{EQOP}). 

If we are in the case where $G_1>2G_0$, we will need to extend first the function $\varphi_{|\mathcal{O}^0}$ 
in the direction of negative $x_2$,  with a finite number of  successive reflections.  That is, if $\varphi_0$ is defined
in $\mathcal{O}^{ \epsilon}$ then we extend $\varphi_0$ to the set $(0,1)\times (-G_0,0)$ with the reflecting along the 
line $x_2=0$.  This produces the function 
$$
\varphi_1 (x_1, x_2) = 
\left\{ 
\begin{array}{ll}
\varphi_0(x_1, x_2) & (x_1, x_2 ) \hbox{  with } 0<x_2<G_\eps(x_1)\\
\varphi_0(x_1, -x_2) & (x_1, x_2 )\hbox{  with } -G_0<x_2\leq 0. 
\end{array}
\right.
$$

We can continue producing these reflections inductively as follows
\begin{equation*}
\varphi_n (x_1, x_2) \!=\! \left\{ 
\begin{array}{ll}
\varphi_{n-1}(x_1,x_2) & (x_1, x_2 )\ \hbox{ with } -(n-1)G_0<x_2\leq G_\eps(x_1) \\
\varphi_{n-1}(x_1,-x_2  - 2(n\!-\!1)G_0) & (x_1, x_2 ) \ \hbox{  with } -nG_0<x_2\leq -(n\!-\!1)G_0\\
\end{array}
\right.
\end{equation*}

Choosing $n$ large enough so that $nG_0>G_1$, constructing $\varphi_n$ and applying to $\varphi_n$ the
procedure by reflection in the $x_2$ direction along the oscillating boundary, we obtain the extension operator $P_\eps$ which
satisfies (\ref{EQOP}). 	
\end{proof}

\begin{remark}\label{rem:stoperator}
1)  This operator preserves periodicity in the $x_1$ variable. That is, if the function $\varphi_\eps$ is periodic in $x_1$, then 
the extended function $P_\eps\varphi_\eps$ is also periodic in $x_1$. 

2) This result can be applied to the case $G_\eps(x)=G(x)$ independent of $\eps$. In particular, we can
apply the extension operator to the basic cell. 
\end{remark}

Now we are in contidion to state our main result.
We consider the general case, that is, the domain $\Omega^\eps$ is given as 
$$
\Omega^\epsilon = \{ (x_1,x_2) \in \R^2 \; | \;  x_1 \in I,  \quad
 0 < x_2 < G_\epsilon(x_1) \}.
$$
where the function $G_\eps(\cdot)$ satisfies hypothesis (H). Recall that we denote by $\Omega=(0,1)\times (0,G_1)$ 
and 
$$\hat \Omega=\Omega \setminus \cup_{i=1}^N\{ (\xi_i,x_2):   G_0<x_2<G_1\},$$
where the points $\xi_0,\xi_1,\ldots, \xi_N$ are given by (H). 

\begin{theorem} \label{GT}
Let $u^\epsilon$ be the  solution of (\ref{P}) with $f^\epsilon \in L^2(\Omega^\epsilon)$ satisfying
$\| f^\epsilon \|_{L^2(\Omega^\epsilon)} \le C$
for some positive constant $C$ independent of $\epsilon > 0$. Assume that the function $\hat f^\eps(x)=\int_0^{G_\eps(x)}f(x,y)dy$ satisfies that $\hat f^\eps\rightharpoonup \hat f$, w-$L^2(0,1)$. 

Then, there exists $\hat u \in H^1(\hat \Omega)$, such that, if $P_\eps$ is the extension operator constructed in Lemma \ref{EOT}, then  
$$
P_{\epsilon} u^\epsilon \rightharpoonup \hat u \quad w - H^1(\hat \Omega)
$$ 
where $\hat u(x_1,x_2)$ depends only on the first variable, that is, $\hat u(x_1,x_2) = u(x_1)$,  and $u$ 
is the unique solution of the Neumann problem
\begin{equation} \label{VFPDL}
\int_{I} \Big\{ r(x)  \, u_x(x) \, \varphi_x(x) 
+ p(x) \, u(x) \, \varphi(x) \Big\} dx = \int_{I}  \, \hat f(x) \, \varphi(x) \, dx
\end{equation}
for all $\varphi \in H^1(I)$, where 
\begin{equation} \label{RPFL}
\begin{gathered}
r(x) =  \frac{1}{L}\int_{Y^*(x)}\Big\{ 1 - \frac{\partial X(x)}{\partial y_1}(y_1,y_2) \Big\} dy_1 dy_2\\ 
p(x) = \frac{|Y^*(x)|}{L} 
\end{gathered}
\end{equation}
and $X(x)$ is the unique solution of problem 
\begin{equation} \label{AUXG}
\left\{
\begin{gathered}
- \Delta X(x)  =  0  \textrm{ in } Y^*(x)  \\
\frac{\partial X(x)}{\partial N}  =  0  \textrm{ on } B_2(x)  \\
\frac{\partial X(x)}{\partial N}  =  N_1 \textrm{ on } B_1(x)  \\
X(x)(0,y_2)  =  X(x)(L,y_2) \textrm{ on } B_0(x)  \\
\int_{Y^*(x)} X(x) \; dy_1 dy_2  =  0  
\end{gathered}
\right.
\end{equation}
in the representative cell $Y^*(x)$ given by
\begin{equation} \label{CELLL}
Y^*(x) = \{ (y_1,y_2) \in \R^2 \; : \; 0< y_1 < L, \quad 0 < y_2 < G(x,y_1) \}. 
\end{equation}
$B_0(x)$ is the lateral boundary, $B_1(x)$ is the upper boundary and $B_2(x)$ is the lower
boundary of $\partial Y^*(x)$ for all $x \in I$.
\end{theorem}

\begin{remark}
i) If the function $ r(x)$ is a continuous function, then, the integral formulation (\ref{VFPDL}) 
is the weak formulation of problem (\ref{GLP-cont}) with $f(x)=\hat f(x)/p(x)$
\par\noindent ii) If initially the function $f^\eps(x,y)=f_0(x)$, then it is not difficult to see that $\hat f^\eps(x)=G_\eps(x)f_0(x)$ and $G_\eps(x)\rightharpoonup \frac{|Y^*(x)|}{L}\equiv p(x)$ as $\eps\to 0$,  and therefore, $\hat f(x)=p(x)f_0(x)$.

\end{remark}

\section{The piecewise periodic case}
\label{piecewise-periodic}

In this section  we find the limit of the sequence $\{ u^\epsilon \}_{\epsilon > 0}$ given by
the Neumann problem (\ref{P}) as $\epsilon$ goes to zero for 
the case where the oscillating boundary is piecewise periodic. 

\begin{figure}[h] 
\centering {\includegraphics[width=7cm,height=6cm]{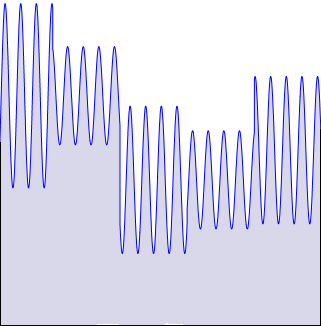}}
\caption{Piecewise periodic domain $\Omega^\eps$}
\label{figura2}
\end{figure}

So let us consider that the family of domains $\Omega^\eps$ satisfies (H) and morever the 
function $G$ is independent of the first variable in each of the domains $(\xi_{i-1},\xi_i)\times \R$. That is, 
there exist $0=\xi_0<\xi_1<\ldots<\xi_{N-1}<\xi_{N}=1$ so that 
the function $G$ from \eqref{def-G} satisfies that  $G(x, y)=G_i(y)$ for $x\in I_i=(\xi_{i-1},\xi_i)$ and $G_i(y+L)=G_i(y)$ for
all $y\in \R$. Moreover, the function $G_i(\cdot)$ is $C^1$ for all $i=1,2,\ldots N$, there exists $0<G_0<G_1$ such that
$0<G_0\leq G_i(\cdot)\leq G_1$ for all $i=1,\ldots, N$;  and the domain is given by
\begin{equation}\label{def-Omega-eps}
\begin{array}{l}
\Omega^\eps=\Big\{ (x,y):  \xi_{i-1}<x<\xi_i, 0<y<G_i(x/\eps), i=1,\ldots, N\Big\}\cup \\ \\
\displaystyle \qquad\qquad \cup_{i=1}^{N-1}\Big\{(\xi_i, y), 0<y<\min\{G_{i-1}(\xi_i/\eps), G_{i}(\xi_i/\eps)\}\Big\}
\end{array}
\end{equation}
(see Figure \ref{figura2}).  Denote also by $\Omega$ the open rectangle $\Omega=I\times (0, G_1)$ and by 
$$\hat \Omega=\Omega \setminus \cup_{i=1}^N\{ (\xi_i,x_2):   G_0<x_2<G_1\}.$$

Observe that for $\Omega^\eps$ and $\hat\Omega$ we have the extension operator $P_\eps$
constructed in Lemma \ref{EOT}.

We can show, 

\begin{theorem} \label{PPCT}
Assume that  $f^\epsilon \in L^2(\Omega^\epsilon)$ satisfies (\ref{ESTF}) and the function $\hat f^\eps(x)=\int_0^{G_\eps(x)}f(x,y)dy$ satisfies that $\hat f^\eps\rightharpoonup \hat f$, w-$L^2(0,1)$. 
Let $u^\epsilon$ be the unique solution of (\ref{P}). 
Then, there exists $\hat u \in H^1(\Omega)$ such that if $P_\epsilon$ is the extension operator constructed in Lemma \ref{EOT},  we have 
$$
P_{\epsilon} u^\epsilon \eto \hat u \quad w - H^1(\hat\Omega), \, s-L^2(\hat\Omega)
$$ 
where $\hat u(x_1,x_2) = u(x_1)$ in $\hat\Omega$ and $u(\cdot)$  is the unique 
weak solution of the Neumann problem
\begin{equation} \label{VFPDL-piecewise0}
\int_{I} \Big\{ r(x)  \; u_x(x) \, \varphi_x(x) 
+ p(x) \, u(x) \, \varphi(x) \Big\} dx = \int_{I}  \hat f(x) \, \varphi(x) \, dx
\end{equation}
for all $\varphi \in H^1(I)$, 
where $p(x)$ and $r(x)$ are piecewise constant functions defined as follows: 
$p(x)=p_i$ for all $x\in (\xi_{i-1},\xi_i)$ where 
\begin{equation}\label{def-p}
p_i=\frac{|Y_i^*|}{L},  \quad i=1,\ldots,N
\end{equation}

where  $Y_i^*$ is the basic cell for $x\in(\xi_{i-1},\xi_i)$, that is
 \begin{equation} \label{CELL}
Y^*_i = \{ (y_1,y_2) \in \R^2 \; : \; 0< y_1 < L, \quad 0 < y_2 < G_i(y_1) \}
\end{equation}
and   $r(x)=r_i$ for all $x\in (\xi_{i-1},\xi_i)$ where
$$
r_i = \frac{1}{L} \int_{Y^*_i} \Big\{ 1 - \frac{\partial X_i}{\partial y_1}(y_1,y_2) \Big\} dy_1 dy_2
$$
and the function $X_i$ is the unique solution of
\begin{equation} \label{AUX}
\left\{
\begin{gathered}
- \Delta X_i  =  0  \textrm{ in } Y_i^*  \\
\frac{\partial X_i}{\partial N}  =  0  \textrm{ on } B^i_2  \\
\frac{\partial X_i}{\partial N}  =  N_1 \textrm{ on } B^i_1  \\
X_i(0,y_2)  =  X_i(L,y_2) \textrm{ on } B^i_0  \\
\int_{Y_i^*} X_i \; dy_1 dy_2  =  0  
\end{gathered}
\right.
\end{equation}

\end{theorem}

\begin{remark} If we define $f_0(x)=\hat f(x)/p(x)$, then problem (\ref{VFPDL-piecewise0}) is equivalent to the following:

\begin{equation} \label{EQF2}
- q_i u_{xx}(x) + u(x) = f_0(x) \quad  x \in (\xi_{i-1},\xi_i)\\
\end{equation}
for $i = 1, ..., N$, where $q_i=r_i/p_i$,  satisfying the following boundary conditions
\begin{equation} \label{BCF2}
\left\{
\begin{gathered}
u_x(\xi_0) = u_x(\xi_N) = 0 \\
q_i \, u_x(\xi_i-) - q_{i+1} \, u_x(\xi_i+) = 0 \quad i = 1, ..., N -1.
\end{gathered}
\right. 
\end{equation} 
Here, $u_x(\xi_i \pm)$ denote the right(left)-hand side limits of $u_x$ at $\xi_i$.

\end{remark}

\begin{proof}  
Let us consider the family of representative cell $Y^*_i$, $i=1,2\ldots,N$,  defined by 
$$
Y^*_i=\{ (y_1,y_2) \in \R^2 \; : \; 0<y_1<L \textrm{ and } 0<y_2<G_i(y_1)\}
$$
and let $\chi_i$ be their characteristic function.
We extend each $\chi_i$ periodically on the variable $y_1 \in \R$ and denote this extension again by $\chi_i$, for $i=1,\ldots, N$.  

If we denote by $\chi^\epsilon_i$ the characteristic function of the set $\Omega^\epsilon_i=\{(x,y):  \xi_{i-1}<x<\xi_{i}, \, 0<y<G_i(x/\eps)=G(x,x/\eps)\}$, 
we easily see that
\begin{equation} \label{chip}
\chi^\epsilon_i(x_1,x_2) 
= \chi_i (\frac{x_1}{\epsilon},x_2),\hbox{ for }(x_1,x_2)  \in \Omega^\epsilon_i.
\end{equation}
Let us also denote by $\Omega^i$ the rectangle 
$\Omega^i=\{(x,y): \xi_{i-1}<x<\xi_i, 0<y<G_1\}$, for $i=1,\ldots,N$. 

Let us also consider the following families of isomorphisms 
$T^{\epsilon}_k : A^{\epsilon}_k \mapsto Y$ given by
\begin{equation} \label{ISO}
T^{\epsilon}_k(x_1,x_2) = (\frac{x_1 - \epsilon k L}{\epsilon},x_2)
\end{equation}
where 
$$
\begin{gathered}
A^{\epsilon}_k = \{ (x_1,x_2) \in \R^2 \; | \;  
\epsilon  k L \leq x_1 < \epsilon L (k+1)  \textrm{ and } 0 < x_2 < G_1 \} \\
Y = (0, L) \times (0, G_1)
\end{gathered}
$$
with $k \in \N$.

Let us consider the following auxiliary problem given by
\begin{equation} \label{AUX0}
\left\{
\begin{gathered}
- \Delta X_i  =  0  \textrm{ in } Y^*_i \\
\frac{\partial X_i}{\partial N}  =  0  \textrm{ on } B_2^i  \\
\frac{\partial X_i}{\partial N} 
= - \frac{G_i'(y_1)}{\sqrt{1+G_i'(y_1)^2}} \textrm{ on } B_1^i  \\
 X_i(0,y_2) = X_i(L,y_2) \textrm{ on } B_0^i\\
  \int_{Y^*_i} X_i \; dy_1 dy_2 = 0  
 \end{gathered}
\right.
\end{equation}
where $B_0^i$ is the lateral boundary, $B_1^i$ is the upper boundary and $B_2^i$ is the lower boundary of $\partial Y^*_i$.

Taking the isomorphism (\ref{ISO}) and the family of extension operators 
$$
P^i \in \mathcal{L}(H^1(Y^*_i),H^1(Y)) \cap \mathcal{L}(L^2(Y^*_i),L^2(Y))
$$ 
defined by Lemma \ref{EOT} with $G_\epsilon(x_1) = G_i(x_1)$
independent of $\epsilon$ and $Y = (0, L) \times (0, G_1)$, see Remark \ref{rem:stoperator},
we define the function
\begin{eqnarray*}
\omega^{\epsilon} (x_1,x_2) & = & x_1 - \epsilon \Big(  P^iX_i \circ T^\epsilon_k (x_1,x_2) \Big) \\
& = & x_1 - \epsilon \Big( P^i X_i (\frac{x_1 - \epsilon L k}{\epsilon},x_2) \Big),  \hbox{ for } (x_1,x_2)\in \Omega^i\cap A^\eps_k, \quad i=1,\ldots, N.
\end{eqnarray*}

Notice that this function is  defined in  $\cup_{i=1}^N\Omega^i$ and it is well defined. For $\eps>0$ fixed and for 
$(x_1,x_2)\in \Omega^i$ for some $i=1,\ldots,N$ then there exists a unique $k\in \N$ such that $(x_1,x_2)\in A^\eps_k$. Observe 
that $\cup_{i=1}^N\Omega^i=\Omega\setminus\cup_{i=1}^{N-1}\{(\xi_i,y): 0<y<G_1\}$ and that $\omega^\eps\in H^1(\cup_{i=1}^N\Omega^i)$.

We introduce now the vector $\eta^\epsilon = (\eta_1^\epsilon,\eta_2^\epsilon)$ defined by 
\begin{equation} \label{VETA}
\eta_r^\epsilon(x_1,x_2) = \frac{\partial \omega^\epsilon}{\partial x_r}(x_1,x_2),  \quad (x_1,x_2)\in \cup_{i=1}^N\Omega^i, \, r=1,2
\end{equation}
Since $\frac{\partial}{\partial x_1} = \frac{1}{\epsilon} \frac{\partial}{\partial y_1}$ 
and $\frac{\partial}{\partial x_2} = \frac{\partial}{\partial y_2}$ we have that
$$
\begin{gathered}
\eta_1^\epsilon(x_1,x_2)= 1 - \frac{\partial X_i}{\partial y_1}(\frac{x_1 - \epsilon k L}{\epsilon},x_2)
= 1 - \frac{\partial X }{\partial y_1}(\frac{x_1}{\epsilon},x_2) := \eta_1(y_1,y_2), \\ 
\eta_2^\epsilon(x_1,x_2)= - \epsilon \frac{\partial X_i }{\partial y_2}(\frac{x_1 - \epsilon k L}{\epsilon},x_2)
= - \epsilon \frac{\partial X }{\partial y_2}(\frac{x_1}{\epsilon},x_2)
:= \eta_2(y_1,y_2)
\end{gathered}
$$
for $(y_1,y_2) = (\frac{x_1 - \epsilon k L}{\epsilon},x_2) \in Y^*_i$ and $(x_1,x_2) \in \Omega_i^\epsilon$, for $i=1,\ldots,N$.

It follows from definition of $X$ that the functions $\eta_1^\eps$ and $\eta_2^\eps$ satisfy
\begin{equation}\label{DIV}
\begin{gathered}
\frac{\partial \eta_1^\epsilon}{\partial {x_1}} 
+ \frac{1}{\epsilon^2} \frac{\partial \eta_2^\epsilon}{\partial{x_2}}   =  0 \textrm{ in }  \Omega^\epsilon_i,\, i=1,\ldots,N \\ 
\eta_1^\epsilon N^\epsilon_1 
+ \frac{1}{\epsilon^2} \eta_2^\epsilon N^\epsilon_2 = 0 \textrm{ on } (x_1, G_i \Big(\frac{x_1}{\epsilon} \Big) ),\, i=1,\ldots,N \\
\eta_1^\epsilon N^\epsilon_1 
+ \frac{1}{\epsilon^2} \eta_2^\epsilon N^\epsilon_2 = 0 \textrm{ on } (x_1, 0) 
\end{gathered}
\end{equation}
where 
$$
\begin{gathered}
N^\epsilon = (N^\epsilon_1, N^\epsilon_2) = \left(
 - \frac{G_i'(\frac{x_1}{\epsilon})}{(\epsilon^2+{G_i'(\frac{x_1}{\epsilon})}^2)^{\frac{1}{2}}}, 
 \frac{\epsilon}{ (\epsilon^2+{G_i'(\frac{x_1}{\epsilon})}^2)^{\frac{1}{2}}} \right) 
 \textrm{ on } (x_1, G_i \Big(\frac{x_1}{\epsilon} \Big) ),\, i=1,\ldots,N \\
N^\epsilon = ( 0 , -1) \textrm{ on } (x_1, 0). 
\end{gathered}
$$

In fact, by (\ref{AUX0}) we have
$$
\frac{\partial \eta_1^\epsilon}{\partial {x_1}} 
+ \frac{1}{\epsilon^2} \frac{\partial \eta_2^\epsilon}{\partial{x_2}}   =  
- \frac{1}{\epsilon} \Big( \frac{\partial^2 X_i}{{\partial y_1}^2}(\frac{x_1}{\epsilon},x_2)
+ \frac{\partial^2 X_i}{{\partial y_2}^2}(\frac{x_1}{\epsilon},x_2) \Big) = 0 \textrm{ in }  \Omega^\epsilon_i, \, i=1,\ldots,N
$$
\begin{eqnarray*}
\eta_1^\epsilon N^\epsilon_1 
+ \frac{1}{\epsilon^2} \eta_2^\epsilon N^\epsilon_2 
& = & - \frac{
 G_i'(\frac{x_1}{\epsilon}) \Big( 1 - \frac{\partial X}{\partial y_1}(\frac{x_1}{\epsilon},x_2) \Big)
+ \frac{\partial X}{\partial y_2}(\frac{x_1}{\epsilon},x_2)  }{\sqrt{ \epsilon^2 + g'(\frac{x_1}{\epsilon})^2} } \\
& = & -  \frac{ G_i'(\frac{x_1}{\epsilon}) 
 + ( \frac{\partial X}{\partial y_1}(\frac{x_1}{\epsilon},x_2), \frac{\partial X}{\partial y_2}(\frac{x_1}{\epsilon},x_2) ) 
 \cdot ( - G_i'(\frac{x_1}{\epsilon}) , 1 ) }
 {\sqrt{ \epsilon^2 + G_i'(\frac{x_1}{\epsilon})^2}}  \\
& = & - \frac{G_i'(\frac{x_1}{\epsilon}) 
- \frac{G_i'(\frac{x_1}{\epsilon})}{\sqrt{ 1 + G_i'(\frac{x_1}{\epsilon})^2}}\sqrt{ 1 + G_i'(\frac{x_1}{\epsilon})^2}}
{\sqrt{ \epsilon^2 + G_i'(\frac{x_1}{\epsilon})^2}}
 = 0 
\end{eqnarray*}
on $(x_1, G_i\Big(\frac{x_1}{\epsilon}\Big) )$ for $i=1,\ldots,N$.  Moreover,
$$
\eta_1^\epsilon N^\epsilon_1 
+ \frac{1}{\epsilon^2} \eta_2^\epsilon N^\epsilon_2 
= ( 1 - \frac{\partial X_i}{\partial y_1}(\frac{x_1}{\epsilon},x_2) ) \cdot 0 
+ \frac{1}{\epsilon}\frac{\partial X_i}{\partial y_2}(\frac{x_1}{\epsilon},x_2) = 0 \textrm{ on } (x_1, 0).
$$
Therefore, multiplying (\ref{DIV}) by a test function $\varphi \in H^1(\Omega)$ with $\varphi=0$ in  neighborhood of the set
$\cup_{i=0}^N \{ (\xi_i,x_2):  0\leq x_2\leq G_1 \}$ and integrating by parts, 
we obtain
\begin{eqnarray*}
0 & = & \int_{\Omega^\epsilon} \varphi \left( \frac{\partial \eta_1^\epsilon}{\partial {x_1}} 
+ \frac{1}{\epsilon^2} \frac{\partial \eta_2^\epsilon}{\partial{x_2}} \right) dx_1 dx_2 \\
& = &  \int_{\partial \Omega^\epsilon}  \varphi \left( \eta_1^\epsilon N^\epsilon_1 + \frac{1}{\epsilon^2} \eta_2^\epsilon N^\epsilon_2 \right) dS 
- \int_{\Omega^\epsilon} \left(\frac{\partial \varphi}{\partial x_1} \eta_1^\epsilon 
+ \frac{1}{\epsilon^2} \frac{\partial \varphi}{\partial x_2} \eta_2^\epsilon \right) dx_1 dx_2 \\ 
& = &  \int_{(x_1, g(\frac{x_1}{\epsilon}) )}  
\varphi \left( \eta_1^\epsilon N^\epsilon_1 + \frac{1}{\epsilon^2} \eta_2^\epsilon N^\epsilon_2 \right) dS 
- \int_{\Omega^\epsilon} \left(\frac{\partial \varphi}{\partial x_1} \eta_1^\epsilon 
+ \frac{1}{\epsilon^2} \frac{\partial \varphi}{\partial x_2} \eta_2^\epsilon \right) dx_1 dx_2 \\
& = & 0 - \int_{\Omega^\eps}  \left(\frac{\partial \varphi}{\partial x_1}  \eta_1^\epsilon  
+ \dfrac{1}{\epsilon^2} \frac{\partial \varphi}{\partial x_2}  \eta_2^\epsilon \right) dx_1 dx_2.
\end{eqnarray*}

That is
\begin{equation}\label{equation-final}
\int_{\Omega^\eps}  \left(\eta_1^\epsilon   \frac{\partial \varphi}{\partial x_1}  
+  \eta_2^\epsilon \dfrac{1}{\epsilon^2} \frac{\partial \varphi}{\partial x_2} \right) dx_1 dx_2=0.
\end{equation}

Let $\phi= \phi(x_1) \in \mathcal{C}^\infty_0(\cup_{i=0}^{N-1} (\xi_i,\xi_{i+1}))$ and considering the test function  $\varphi=\phi \; \omega^\epsilon$ in (\ref{VFP})  and in (\ref{equation-final}), we obtain

\begin{eqnarray}
& & \int_{\Omega^\eps} f^\eps (\phi \; \omega^\epsilon) dx_1 dx_2 \nonumber \\
& = & 
\int_{\Omega^\eps} \Big\{ {\frac{\partial u^\epsilon}{\partial x_1}} \frac{\partial}{\partial x_1} (\phi \; \omega^\epsilon)  + \frac{1}{\epsilon^2} {\frac{\partial u^\epsilon}{\partial x_2}} \frac{\partial}{\partial x_2} (\phi \; \omega^\epsilon) 
+  {u^\epsilon} (\phi \; \omega^\epsilon) \Big\} dx_1 dx_2\nonumber \\
& = & 
\int_{\Omega^\eps} \Big\{ {\frac{\partial u^\epsilon}{\partial x_1}} \frac{\partial}{\partial x_1} (\phi \; \omega^\epsilon)  + 
\frac{1}{\epsilon^2} {\frac{\partial u^\epsilon}{\partial x_2}} \frac{\partial}{\partial x_2} (\phi \; \omega^\epsilon) 
+  {u^\epsilon} (\phi \; \omega^\epsilon) \Big\} dx_1 dx_2\nonumber \\
& & -  \int_{\Omega^\eps} \Big\{ {\eta_{1}^\epsilon} \frac{\partial }{\partial x_1} (\phi \;  \omega^\epsilon)  
+ \frac{1}{\epsilon^2}  {\eta_{2}^\epsilon}  \frac{\partial}{\partial x_2} (\phi \, \omega^\epsilon) \Big\} dx_1 dx_2 \nonumber \\
& = & \int_{\Omega^\eps} \Big\{ {\frac{\partial u^\epsilon}{\partial x_1}} 
\frac{\partial \phi}{\partial x_1}   \omega^\epsilon 
+ {\frac{\partial u^\epsilon}{\partial x_1}} \frac{\partial \omega^\epsilon}{\partial x_1}  \phi 
+ \frac{1}{\epsilon^2} {\frac{\partial u^\epsilon}{\partial x_2}} \frac{\partial \omega^\epsilon}{\partial x_2}  \phi 
+ {u^\epsilon} \phi  \omega^\epsilon \nonumber \\
& & - {\eta_{1}^\epsilon}  \frac{\partial \phi}{\partial x_1}  u^\epsilon
-  {\eta_{1}^\epsilon} \frac{\partial  u^\epsilon}{\partial x_1}  \phi
- \frac{1}{\epsilon^2} {\eta_{2}^\epsilon} \frac{\partial  u^\epsilon}{\partial x_2}  \phi \Big\} dx_1 dx_2.
\label{EQ!}
\end{eqnarray}

Using that $\eta_i^\epsilon=\frac{\partial \omega^\epsilon}{\partial x_i} $ we cancel the appropriate terms and obtain 
\begin{equation} \label{EQ08}
\begin{array}{l}
\displaystyle
\int_{\Omega^\eps} \Big\{{\frac{\partial u^\epsilon}{\partial x_1}} \frac{\partial \phi}{\partial x_1} \, \omega^\epsilon 
-  {\eta_1^\epsilon}  \frac{\partial \phi}{\partial x_1} \,   u^\epsilon 
+  {u^\epsilon} \phi \, \omega^\epsilon
\Big\}  dx_1 dx_2 \\ \\
\displaystyle \qquad\qquad = \int_{\Omega^\eps} f^\eps \phi \, \omega^\epsilon dx_1 dx_2,\quad \forall \phi\in  \mathcal{C}^\infty_0(\cup_{i=0}^{N-1} (\xi_i,\xi_{i+1})).
\end{array}
\end{equation}

On the other hand, we have obtained before the weak formulation of problem (\ref{P}), that is 
\begin{equation} \label{VFP0-bis1}
\int_{\Omega^\epsilon} \Big\{ \frac{\partial u^\epsilon}{\partial x_1} \frac{\partial \varphi}{\partial x_1} 
+ \frac{1}{\epsilon^2} \frac{\partial u^\epsilon}{\partial x_2} \frac{\partial \varphi}{\partial x_2}
+ u^\epsilon \varphi \Big\} dx_1 dx_2 = \int_{\Omega^\epsilon} f^\eps \varphi dx_1 dx_2, 
\,\, \forall \varphi \in H^1(\Omega^\epsilon).
\end{equation}

Now we need to pass to the limit in (\ref{EQ08}) and (\ref{VFP0-bis1}). In order to accomplish this we need
to write both expressions as integrals in the same domain. For this, we will use the extension $P_\eps$ constructed
in Lemma  \ref{EOT}, the standard extension by zero, that we denote by $\widetilde{}\, $, and the characteristic function $\chi^{\eps}$ of $\Omega^{\eps}$ as follows:

\begin{equation} \label{EQ08-bis}
\begin{array}{l}
\displaystyle
\int_{\Omega} \Big\{\widetilde{{\frac{\partial u^\epsilon}{\partial x_1}}} \frac{\partial \phi}{\partial x_1} \, \omega^\epsilon 
-  \widetilde{\eta_1^\epsilon}  \frac{\partial \phi}{\partial x_1} \,  P_\eps u^\epsilon 
+  \chi^{\eps} P_\eps({u^\epsilon}) \phi \, \omega^\epsilon
\Big\}  dx_1 dx_2\qquad\qquad\qquad  \\ \\
\displaystyle \qquad \qquad\qquad = \int_{\Omega} \chi^{ \eps} f^\eps \phi \, \omega^\epsilon dx_1 dx_2,\qquad \forall \phi\in \mathcal{C}^\infty_0(\cup_{i=0}^{N-1} (\xi_i,\xi_{i+1}))
\end{array}
\end{equation}

\begin{equation} \label{VFP0-bis}
\int_{\Omega} \Big\{ \widetilde{\frac{\partial u^\epsilon}{\partial x_1}} \frac{\partial \varphi}{\partial x_1} 
+ \chi^{\eps} P_\eps u^\epsilon \varphi \Big\} dx_1 dx_2 = \int_{\Omega}\chi^{\eps} f^\eps \varphi dx_1 dx_2, 
\,\, \forall \varphi \in H^1(0,1).
\end{equation}
Observe that in this last equality we have taken $\varphi \in H^1(0,1)$ and the term including partial derivatives with
respect to $x_2$ do not appear. 

We want to pass to the limit in the expressions above, (\ref{EQ08-bis}) and (\ref{VFP0-bis}) . In order to accomplish this, we pass to the limit in the different
functions that form the integrands. 

\par\bigskip\noindent {\bf  (a). Limit of $\chi_\eps$}. 

From (\ref{chip}), we have for $i=1,\ldots,N$, 
\begin{equation} \label{chi0}
\chi^\epsilon_i( \cdot, x_2) \stackrel{\eps\to 0}{\rightharpoonup} \theta_i(x_2) := \frac{1}{L} \int_0^L \chi_i(s,x_2) ds 
\quad w^*-L^\infty((\xi_{i-1},\xi_i), \quad\forall x_2\in (0,G_1).
\end{equation}

Observe that the limit $\theta_i$ does not dependent on the variable $x_1\in (\xi_{i-1},\xi_i)$, although it depends on $i=1,\ldots,N$.  Moreover, 
we can get the area of the open set $Y^*_i$ with the formula
\begin{equation} \label{ITHETA}
L \int_0^{G_1} \theta_i(x_2) dx_2 = |Y^*_i|.
\end{equation}
Also, from (\ref{chi0}) we have that 
$$
H^\epsilon_i(x_2) = \int_{\xi_{i-1}}^{\xi_i}  \varphi(x_1,x_2) \, \Big\{ \chi^\epsilon_i(x_1,x_2) - \theta_i(x_2) \Big\} \, dx_1 
\to 0 \textrm{ as } \epsilon \to 0 
$$ 
a.e. $x_2 \in (0, G_1)$ and for all $\varphi \in L^1(\Omega)$.
So, due to
$$
\begin{gathered}
\int_{\Omega_i} \varphi(x_1,x_2) \,\Big\{ \chi^\epsilon_i(x_1,x_2) - \theta_i(x_2) \Big\} \, dx_1 dx_2 
= \int_0^{G1} H^\epsilon_i(x_2) dx_2 \\ 
\textrm{ and }  |H^\epsilon_i(x_2)| \le \int_{\xi_{I-1}}^{\xi_i}| \varphi(x_1,x_2)| dx_1,
\end{gathered}
$$
we can get by Lebesgue's Dominated Convergence Theorem that
\begin{equation} \label{CHIAIM0}
\chi^\epsilon_i\stackrel{\eps\to 0}{ \rightharpoonup} \theta \quad w^*-L^\infty(\Omega)
\end{equation}
where $\theta(x_1,x_2)=\theta_i(x_2)$ for $x_1\in (\xi_{i-1},\xi_i)$, $i=1,2,\ldots,N$.

\par\bigskip\noindent {\bf (b). Limit in the tilde functions}

Since $\|f^{\epsilon}\|_{L^{2}(\Omega)}$ is uniformly bounded, we get from (\ref{priori})  that there exists $M$ independent of $\epsilon$ such that
\begin{equation} \label{estimate0}
\begin{gathered}
\| \widetilde{u^\epsilon} \|_{L^2(\Omega)}, \Big\| \widetilde{\frac{\partial u^\epsilon}{\partial x_1}} \Big\|_{L^2(\Omega)} \textrm{ and }
\frac{1}{\epsilon} \Big\| \widetilde{\frac{\partial u^\epsilon}{\partial x_2}} \Big\|_{L^2(\Omega)} \le M 
\textrm{ for all } \epsilon > 0.
\end{gathered}
\end{equation}
Then, we can extract a subsequence, still denoted by $\widetilde{u^\epsilon}$, $\widetilde{\frac{\partial u^\epsilon}{\partial x_1}}$
and $\widetilde{\frac{\partial u^\epsilon}{\partial x_2}}$, such that
\begin{equation} \label{WC0}
\begin{gathered}
\widetilde{u^\epsilon} \rightharpoonup u^* \quad w-L^2(\Omega) \\
\widetilde{\frac{\partial u^\epsilon}{\partial x_1}} \rightharpoonup \xi^* \quad w-L^2(\Omega) \textrm{ and }\\
\widetilde{\frac{\partial u^\epsilon}{\partial x_2}} \rightarrow 0 \quad s-L^2(\Omega) \\
\end{gathered}
\end{equation}
as $\epsilon \to 0$ for some  $u^*$ and $\xi^* \in L^2(\Omega)$.

Moreover, since $\|f^{ \eps}\|_{L^2(\Omega^{ \eps})}\leq C$ independent of $\eps$, we have $\|\tilde f^{\eps}\|_{L^2(\Omega)}\leq C$ and
therefore, the function $\hat f^\eps$ defined by 
\begin{equation}\label{def-hat-f}
\hat f^\eps (x_1)\equiv \int_0^{G_1} \tilde f^\eps(x_1,x_2)dx_2  
\end{equation}
satisfies that $\hat f^\eps\in L^2(0,1)$.  Hence, via subsequences, we have the existence of a function $\hat f=\hat f(x_1)\in L^2(0,1)$ such that 
\begin{equation}\label{convergence-f}
\hat f^\eps \rightharpoonup \hat f\qquad w-L^2(0,1).
\end{equation}

\par\bigskip\noindent {\bf (c).  Limit in the extended functions}

Using the a priori estimate (\ref{priori}), the fact that $u^{ \eps} \in H^1(\Omega^{ \eps})$ and using the results from Lemma \ref{EOT} on
the extension operator $P_\epsilon$ we get that 
\begin{equation} \label{estimateP}
\begin{gathered}
\| P_{\epsilon} u^\epsilon \|_{L^2(\Omega)}, \Big\| \frac{\partial P_{\epsilon} u^\epsilon}{\partial x_1} \Big\|_{L^2(\Omega)} \textrm{ and }
\frac{1}{\epsilon} \Big\| \frac{\partial P_{\epsilon} u^\epsilon}{\partial x_2} \Big\|_{L^2(\Omega)} \le \tilde{M} 
\textrm{ for all } \epsilon > 0
\end{gathered}
\end{equation}
where $\tilde{M}$ is a positive constant independent of $\epsilon$ given by
estimate (\ref{estimate0}) and Lemma \ref{EOT}.
Then, we can extract a subsequence, still denoted by $P_{\epsilon} u^\epsilon$ and a function $ u_0\in H^1( \Omega)$,  such that
\begin{equation} \label{LEO}
\begin{gathered}
P_{\epsilon} u^{\epsilon} \rightharpoonup  u_0 \quad w-H^1( \hat \Omega) \\
P_{\epsilon} u^{\epsilon} \rightarrow  u_0 \quad s-L^2(\hat \Omega) \\
\textrm{ and }
\frac{\partial P_{\epsilon} u^\epsilon}{\partial x_2} \rightarrow 0 \quad s-L^2(\hat \Omega).
\end{gathered}
\end{equation}

A consequence of the limits (\ref{LEO}) is that $ u_0(x_1,x_2)$ does not depend on the variable $x_2$.  
More precisely,
\begin{equation} \label{U0}
\frac{\partial  u_0}{\partial x_2}(x_1,x_2) = 0 \textrm{ a.e. } \Omega.
\end{equation}
In fact, for $i=1,\ldots,N$ and for all $\varphi \in C^\infty_0(\Omega_i)$, we have by (\ref{LEO}) that
\begin{eqnarray*}
\int_{\Omega_i}  u_0 \, \frac{\partial \varphi}{\partial x_2} \, dx_1 dx_2
& = & \lim_{\epsilon \to 0} \int_{\Omega_i} P_{\epsilon} u^{\epsilon} \, \frac{\partial \varphi}{\partial x_2} \, dx_1 dx_2 \\
& = & - \lim_{\epsilon \to 0} \int_{\Omega_i} \frac{\partial P_{\epsilon} u^{\epsilon}}{\partial x_2} \, \varphi \, dx_1 dx_2
= 0
\end{eqnarray*}
which implies that $ u_0(x_1,x_2)$ does not depend on $x_2$. Morever, since the rectangle $I\times (0,G_0)\subset \Omega^\eps$ for all $\eps$ and $u^\eps\in H^1(I\times (0,G_0))$ we have from (\ref{LEO}) that $u^\eps \rightharpoonup  u_0$ 
$w-H^1(I\times (0,G_0))$ and therefore $u_0\in H^1(0,1)$. 

Also, we note that $\widetilde{u^\epsilon} = \chi^\epsilon P_{\epsilon} u^\epsilon \textrm{ a.e. } \Omega$.
Thus, it follows from (\ref{CHIAIM0}), (\ref{WC0}) and (\ref{LEO}) that we have the following relationship between $u^*$ and $u_0$
\begin{equation} \label{U0U*0}
u^*(x_1,x_2) = \theta_i(x_2) \, u_0(x_1) \quad \textrm{ a.e. } (x_1,x_2)\in \Omega_i,\quad i=1,\ldots,N.
\end{equation}

\par\bigskip \noindent {\bf (d). Limit in $\omega_\eps$}. 

With the definition of $\omega_\eps$, we have for all $i=1,\ldots,N$, 

$$
\int_{A^\epsilon_k\cap \Omega^i} |\omega^\epsilon - x_1|^2 dx_1 dx_2 
= \int_{Y_i} \epsilon^3 |(P X_i)(y_1,y_2)|^2 dy_1 dy_2 \le \int_{Y^*_i} C \epsilon^3 |X_i(y_1,y_2)|^2 dy_1 dy_2
$$
and so, 
$$
\begin{gathered}
\int_{\Omega_i} |\omega^\epsilon - x_1|^2 dx_1 dx_2 
\approx  \sum_{k=1}^{\frac{C}{\epsilon L}}  \int_{Y^*_i} C \epsilon^3 |X_i(y_1,y_2)|^2 dy_1 dy_2 \\
\approx  \epsilon^2 \int_{Y^*_i} C |X_i(y_1,y_2)|^2 dy_1 dy_2 \rightarrow 0 \textrm{ as } \epsilon \to 0.
\end{gathered}
$$
Similarly,
\begin{eqnarray*}
\int_{A^\epsilon_k\cap \Omega_i} \Big|\frac{\partial }{\partial x_1} \left( \omega^\epsilon - x_1 \right) \Big|^2 dx_1 dx_2 
& = & \int_{Y_I}  \Big| \frac{\partial (P X_I)}{\partial y_1} (y_1,y_2) \Big|^2 \, \epsilon \, dy_1 dy_2 \\
& \le & \epsilon \int_{Y^*_I} C \Big|\frac{\partial X_i}{\partial y_1}(y_1,y_2)\Big|^2 dy_1 dy_2 
\end{eqnarray*}
and
\begin{eqnarray*}
\int_{A^\epsilon_k\cap \Omega_i} \Big|\frac{\partial }{\partial x_2} \left( \omega^\epsilon - x_1 \right) \Big|^2 dx_1 dx_2 
& = & \int_{Y_i}  \Big|\epsilon \, \frac{\partial (P X_i)}{\partial y_2} (y_1,y_2)\Big|^2 \, \epsilon dy_1 dy_2 \\
& \le & \epsilon^3 \int_{Y^*_i}  C \Big|\frac{\partial X_i}{\partial y_2}(y_1,y_2)\Big|^2 dy_1 dy_2.
\end{eqnarray*}
Also, we have
$$
\begin{gathered}
\int_{\Omega_i} \Big|\frac{\partial }{\partial x_1} \left( \omega^\epsilon - x_1 \right) \Big|^2 dx_1 dx_2 
\approx  \sum_{k=1}^{\frac{C}{\epsilon L}}  \int_{Y^*_i} C \epsilon |\frac{\partial X_i}{\partial y_1}(y_1,y_2)|^2 dy_1 dy_2 \\
\approx  \int_{Y^*_i} \tilde{C} \Big|\frac{\partial X_i}{\partial y_1}(y_1,y_2)\Big|^2 dy_1 dy_2 
\end{gathered}
$$
for all $\epsilon > 0$ and
$$
\begin{gathered}
\int_{\Omega_i} \Big|\frac{\partial }{\partial x_2} \left( \omega^\epsilon - x_1 \right) \Big|^2 dx_1 dx_2 
\le  \epsilon^2 \int_{Y^*_i} \tilde{C} \Big|\frac{\partial X_i}{\partial y_2}(y_1,y_2)\Big|^2 dy_1 dy_2 \to 0 \textrm{ as } \epsilon \to 0.
\end{gathered}
$$

Then, we can conclude
\begin{equation} \label{OMEGAL}
\omega^\epsilon \to x_1 \quad s-L^2(\Omega) \textrm{ and } w-H^1(\Omega_i), \quad i=1,\ldots,N,
\end{equation}
and
$$
\frac{\partial \omega^\epsilon}{\partial x_2} \to 0 \quad s-L^2(\Omega).
$$

\par\bigskip\noindent{\bf (e). Limit of $\eta_1^\eps$}

Let $\widetilde{\eta}^\epsilon = \eta^\epsilon \chi^\epsilon$ be the extension by zero
of the vector $\eta^\epsilon$ to the whole $\Omega$.
We can obtain by the Average Theorem that
\begin{equation} \label{ETAD}
\widetilde{\eta}_1^\epsilon(x_1,x_2) \rightharpoonup \frac{1}{L} \int_0^L  
\Big( 1 - \frac{\partial X_i}{\partial y_1}  (s,x_2) \Big) \chi_i(s,x_2)ds :=  \hat r_i(x_2)
\quad w^*-L^\infty(\xi_{i-1},\xi_i)
\end{equation}
where $\chi_i$ is the characteristic function of $Y^*_i$. 

Hence, arguing as (\ref{CHIAIM0}) we can prove
\begin{equation} \label{ETA}
\widetilde{\eta}_1^\epsilon \rightharpoonup \hat r \quad w^*-L^\infty(\Omega).
\end{equation}
where $\hat r(x_1,x_2)\equiv \hat r_i(x_2)$, $(x_1,x_2)\in \Omega_ i$, $i=1,\ldots,N$.

\par\bigskip\bigskip
Now, by the convergences shown in (a)-(e) above, 
we can pass to the limit in (\ref{EQ08-bis}) and in (\ref{VFP0-bis}).  
We obtain, for all $\phi\in \mathcal{C}^\infty_0(\cup_{i=0}^{N-1} (\xi_i,\xi_{i+1}))$, 
$$
\int_{\Omega} \Big\{ \xi^* \frac{\partial \phi}{\partial x_1} \, x_1 -  \hat  r \frac{\partial \phi}{\partial x_1} \, u_0 +  \theta u_0 \phi \, x_1
\Big\}  dx_1 dx_2 = \int_0^1\hat f(x_1) \phi (x_1)\, x_1 dx_1 .$$
Observe that $\xi^* \frac{\partial}{\partial x_1}(\phi \, x_1) = \xi^* x_1 \frac{\partial \phi}{\partial x_1} + \xi^*  \phi$.
Consequently,  we have
\begin{equation}\label{EQ10}
\int_{\Omega} \Big\{ \xi^* \frac{\partial }{\partial x_1}(\phi \, x_1) - \phi \, \xi^* 
-  \hat r \frac{\partial \phi}{\partial x_1} \, u_0 +  \theta u_0 \phi \, x_1 \Big\} dx_1 dx_2 =  \int_0^1\hat f(x_1) \phi (x_1)\, x_1 dx_1
\end{equation}
for all $\phi\in  \mathcal{C}^\infty_0(\cup_{i=0}^{N-1} (\xi_i,\xi_{i+1}))$. From (\ref{VFP0-bis}), we get 
\begin{equation} \label{EQ100-a}
\int_{\Omega} \Big\{ \xi^* \frac{\partial \varphi}{\partial x_1} +  \theta_i u_0 \varphi \Big\} dx_1 dx_2 
= \int_0^1\hat f \varphi dx_1, \qquad \forall \varphi\in H^1(0,1).
\end{equation}
In particular, via iterated integration and \eqref{EQ100-a}, we get 
\begin{equation} \label{EQ100-a-iterated}
\begin{array}{l}
\displaystyle
\sum_{i=1}^N \int_{\xi_{i-1}}^{\xi_i}\Big\{ \big(\int_0^{G_1}\xi^*(x_1,x_2)dx_2 \big)\frac{\partial \varphi(x_1)}{\partial x_1} + \frac{|Y^*_i|}{L} u_0(x_1)\varphi (x_1)\Big\} dx_1 \\
\displaystyle
\qquad\qquad= \int_0^1\hat f(x_1) \varphi(x_1) dx_1, \qquad \forall \varphi\in H^1(0,1).
\end{array}
\end{equation}

Taking $\varphi=\phi x_1$ in (\ref{EQ100-a}), we get  
\begin{equation} \label{EQ100}
\int_{\Omega} \Big\{ \xi^* \frac{\partial }{\partial x_1}(\phi \, x_1) +  \theta u_0 \phi \, x_1 \Big\} dx_1 dx_2 
= \int_0^1{\hat f \phi x_{1}} dx_1.
\end{equation}
Hence, it follows from (\ref{EQ10}) and (\ref{EQ100}) that, for all $\phi \in\mathcal{C}^\infty_0(\cup_{i=0}^{N-1} (\xi_i,\xi_{i+1}))$,
\begin{equation} \label{EQ11}
0 = \int_{\Omega} \Big\{ \phi \, \xi^* +  \hat r \frac{\partial \phi}{\partial x_1} \,u_0 \Big\} dx_1 dx_2
 =  \int_{\Omega} \Big\{ \phi \, \xi^* - \hat r \phi \frac{\partial u_0}{\partial x_1} \Big\} dx_1 dx_2
\end{equation}
where we have performed an integration by parts to obtain the last integral.  Observe that this integration by parts can
be performed since $\phi \in\mathcal{C}^\infty_0(\cup_{i=0}^{N-1} (\xi_i,\xi_{i+1}))$ and $\hat r$ does not depend on $x_1$ 
 in each of the domains
$\Omega_i$. Hence, if we define
$$
r_i  \equiv \int_{0}^{G_1} \hat r_i(s) ds = \frac{1}{L} \int_{Y^*_i} \Big\{ 1 - \frac{\partial X_i}{\partial y_1}(y_1,y_2) \Big\} \, dy_1 dy_2, 
\quad i=1,\ldots, N
$$
and we denote by 
$$r(x_1)=r_i,\quad\hbox{ for }x_1 \in (\xi_{i-1},\xi_i),\quad i=1,\ldots,N$$ 
and performing an iterated integration in (\ref{EQ11})  we get
$$\int_0^1\phi(x_1)\Big(\int_0^{G_1}\xi^*(x_1,x_2)dx_2-r(x_1) \frac{\partial u_0(x_1)}{\partial x_1} \Big)dx_1=0,\quad \forall \phi\in  \mathcal{C}^\infty_0(\cup_{i=0}^{N-1} (\xi_i,\xi_{i+1}))$$
which implies that
\begin{equation}\label{xi=q}
\int_0^{G_1}\xi^*(x_1,x_2)dx_2=r(x_1) \frac{\partial u_0(x_1)}{\partial x_1},\quad\hbox{ a.e. } x_1\in (0,1).
\end{equation}

Plugging this last equality in (\ref{EQ100-a-iterated}) we get 
\begin{equation} \label{EQ13-bis}
\sum_{i=1}^N\int_{\xi_{i-1}}^{\xi_i} r_i \frac{\partial u_0}{\partial x_1}\frac{\partial \varphi}{\partial x_1} 
+ \frac{|Y^*_i|}{L} \, u_0 \, \varphi dx_1 = \int_0^1  \hat f \, \varphi \, dx_1,\quad \forall \varphi\in H^1(0,1).
\end{equation}

\end{proof}

\section{A domain dependence result} \label{BS}

In this section we are going to analyze how the solutions of \eqref{P} depend on the domain $\Omega^\eps$ and more
exactly on the function $G_\eps$.  As a matter of fact we will show a continuous dependence 
result with respect to the functions $G_\eps$. 

More precisely, assume $G_\eps$ and $\hat G_\eps$ are piecewise continuous functions satisfying (\ref{HG})
and consider the associated oscillating domains $\Omega^\epsilon$ and $\hat \Omega^\epsilon$ given by
$$
\begin{gathered} 
\Omega^\epsilon = \{ (x_1,x_2) \in \R^2 \; | \;  x_1 \in I, \quad 0 < x_2 < G_\epsilon(x_1) \} \\
\hat \Omega^\epsilon = \{ (x_1,x_2) \in \R^2 \; | \;  x_1 \in I, \quad 0 < x_2 < \hat G_\epsilon(x_1) \}.
\end{gathered}
$$

Let $u^\epsilon$ and $\hat u^\epsilon$ be the solutions of the problem (\ref{P}) in the oscillating domains 
$\Omega^\epsilon$ and $\hat \Omega^\epsilon$ respectively with $f^\epsilon \in L^2(\R^2)$.
Then we have the following result:

\begin{theorem}  \label{BPT}
There exists a positive real function $\rho:[0,\infty) \mapsto [0,\infty)$ such that
\begin{equation}\label{main-inequality}
\|u^\epsilon-\hat u^\epsilon\|^2_{H^1_\epsilon(\Omega^\epsilon \cap \hat \Omega^\epsilon)} 
+ \|u^\epsilon\|^2_{H^1_\epsilon(\Omega^\epsilon \setminus \hat \Omega^\epsilon)}
+ \|\hat u^\epsilon\|^2_{H^1_\epsilon(\hat\Omega^\epsilon \setminus \Omega^\epsilon)} \leq \rho(\delta)
\end{equation}
with $\rho(\delta)\to 0$ as $\delta\to 0$ uniformly for all 
\begin{itemize}
\item $\epsilon > 0$;
\item piecewise $C^1$ functions $G_\eps$ and $\hat G_\eps$ with $0<G_0\leq G_\eps(\cdot),\hat G_\eps (\cdot)\leq G_1$ and 
$$
\|G_\eps-\hat G_\eps\|_{L^\infty(0,1)} \leq \delta;
$$
\item $f^\epsilon\in L^2(\R^2)$, $\|f^\epsilon\|_{L^2(\R^2)}\leq 1$.
\end{itemize}
\end{theorem}

\begin{remark} The important part of this result is that the function $\rho(\delta)$ does not depend on $\eps$. Only
depends on $G_0$ and $G_1$.

\end{remark}

To prove this theorem, we use the fact that $u^\epsilon$ and $\hat u^\epsilon$ are minimizers of the quadratic forms
\begin{equation} \label{FUNCAUX}
\begin{gathered}
V_\epsilon(\varphi) =
\frac{1}{2} \int_{\Omega^\epsilon} \Big\{ \frac{\partial \varphi}{\partial x_1}^2 
+ \frac{1}{\epsilon^2} \frac{\partial \varphi}{\partial x_2}^2 
+ \varphi^2 \Big\} dx_1 dx_2 - \int_{\Omega^\epsilon} f^\epsilon \varphi dx_1 dx_2 \\
\hat V_\epsilon(\hat \varphi) =
\frac{1}{2} \int_{\hat \Omega^\epsilon} \Big\{ \frac{\partial \hat \varphi}{\partial x_1}^2 
+ \frac{1}{\epsilon^2} \frac{\partial \hat \varphi}{\partial x_2}^2 
+ \hat \varphi^2 \Big\} dx_1 dx_2 - \int_{\hat \Omega^\epsilon} f^\epsilon \hat \varphi dx_1 dx_2.
\end{gathered}
\end{equation} 
That is,  we have
$$
\begin{gathered}
V_\epsilon(u^\epsilon) = \min_{\varphi \in H^1(\Omega^\epsilon)} V_\epsilon(\varphi) \\
\hat V_\epsilon(\hat u^\epsilon) = \min_{\hat \varphi \in H^1(\hat \Omega^\epsilon)}\hat V_\epsilon(\hat \varphi).
\end{gathered}
$$

In order to prove Theorem \ref{BPT} we will need to consider the minimizer of the functionals $V_\eps$, $\hat V_\eps$,
 and plug them in the other functional. For this, we need to transform the function $u^\eps$ into 
a function defined in $\hat \Omega^\eps$ and the function $\hat u^\eps$ into a function defined in $\Omega^\eps$. 
One possibility is to use some kind of extension operator as the one we have constructed in Section \ref{PRE}. But the problem with 
this approach is  that the norm of the extension operators will depend on the derivatives of the function $G_\eps$ and
$\hat G_\eps$ and therefore it will be very unlikely to prove a results that will depend only on the $L^\infty$ norm of
$G_\eps-\hat G_\eps$. 

In order to transform the function $u^\eps$ (resp. $\hat u^\eps$) into a function defined in $\hat \Omega^\eps$ (resp. $\Omega^\eps$), 
we construct the following operators:
\begin{equation} \label{PETA}
\begin{gathered}
P_{1+\eta} : H^1(U) \mapsto H^1(U(1+\eta)) \\
(P_{1+\eta}\varphi )(x_1, x_2) = \varphi \left( x_1,\frac{x_2}{1+\eta} \right) \quad (x_1,x_2) \in U
\end{gathered}
\end{equation} 
where 
$$
U(1+\eta) = \{ (x_1, (1 + \eta) x_2) \in \R^2 \, | \, (x_1, x_2) \in U \}
$$
and $U \subset \R^2$ is an arbitrary open set.

We also consider the following norm in $H^1(U)$
\begin{equation} \label{NPETA}
\|w\|_{H^1_{\epsilon,1+\eta}(U)}^2 = \frac{1}{1+\eta} \|w\|_{L^2(U)}^2 + \frac{1}{1+\eta} \|w_{x_1}\|_{L^2(U)}^2
+ \frac{1+\eta}{\epsilon^2} \|w_{x_2}\|_{L^2(U)}^2,
\end{equation}
and we can easily see that 
\begin{equation} \label{INPETA}
\|\varphi\|_{H^1_\epsilon(U)}^2=\|P_{1+\eta} \varphi\|_{H^1_{\epsilon,1+\eta}(U(1+\eta))}^2
\end{equation}
and
\begin{equation} \label{INPETA2}
\frac{1}{1+\eta}\|\varphi\|_{H^1_\epsilon(U)}^2\leq \|\varphi\|_{H^1_{\epsilon,1+\eta}(U)}^2\leq (1+\eta)
\|\varphi\|_{H^1_\epsilon(U)}^2.
\end{equation}

We have the following preliminary result about the behavior of the solutions near of the oscillating boundary.

\begin{lemma} \label{BSL}
Let $u^\epsilon$ be the solution of the problem (\ref{P}) and let $P_{1+\eta}$ be the operator given by (\ref{PETA}).

Then exists a positive constant $C = C(G_1,\| f^\epsilon \|_{L^2})$ independent of $\epsilon \in (0,1)$ such that 
$$
\|u^\epsilon\|_{H^1_\epsilon(\Omega^\epsilon \setminus \Omega^\epsilon(\frac{1}{1+\eta}))}^2 
+ \|P_{1+\eta}u^\epsilon\|_{H^1_\epsilon(\Omega^\epsilon(1+\eta) \setminus \Omega^\epsilon)}^2 +\|u^\epsilon-P_{1+\eta}u^\epsilon\|^2_{H^1_\epsilon(\Omega^\epsilon)} 
\leq C \sqrt{\eta}
$$ 
for all $\eta > 0$.
\end{lemma}
\begin{proof}
Since $\eta > 0$, we have $\Omega^\epsilon(\frac{1}{1+\eta}) \subset \Omega^\epsilon$. 
So, we obtain,
\begin{eqnarray}
V_\epsilon(u^\epsilon)
& = & \frac{1}{2}\|u^\epsilon\|_{H^1_\epsilon(\Omega^\epsilon)}^2-\int_{\Omega^\epsilon}f^\epsilon u^\epsilon \, dx_1 dx_2 \nonumber \\
& = & \frac{1}{2}\|u^\epsilon\|_{H^1_\epsilon(\Omega^\epsilon\setminus \Omega^\epsilon(\frac{1}{1+\eta}))}^2+
\frac{1}{2}\|u^\epsilon\|_{H^1_\epsilon( \Omega^\epsilon(\frac{1}{1+\eta}))}^2-\int_{\Omega^\epsilon}f^\epsilon u^\epsilon \, dx_1 dx_2 \nonumber \\
& = & \frac{1}{2}\|u^\epsilon\|_{H^1_\epsilon(\Omega^\epsilon\setminus \Omega^\epsilon(\frac{1}{1+\eta}))}^2+
\frac{1}{2}\|P_{1+\eta}u^\epsilon\|_{H^1_{\epsilon,1+\eta}( \Omega^\epsilon)}^2-\int_{\Omega^\epsilon}f^\epsilon u^\epsilon \, dx_1 dx_2 \nonumber \\
& \geq & \frac{1}{2}\|u^\epsilon\|_{H^1_\epsilon(\Omega^\epsilon\setminus \Omega^\epsilon(\frac{1}{1+\eta}))}^2+
\frac{1}{2(1+\eta)}\|P_{1+\eta}u^\epsilon\|_{H^1_{\epsilon}( \Omega^\epsilon)}^2-\int_{\Omega^\epsilon}f^\epsilon u^\epsilon\, dx_1 dx_2 \nonumber \\
\label{EQ1}
\end{eqnarray}
where we have used (\ref{INPETA}), (\ref{INPETA2}). Moreover, 
$$
\begin{array}{l}
\|P_{1+\eta}u^\epsilon\|_{H^1_{\epsilon}( \Omega^\epsilon)}^2 =\|P_{1+\eta}u^\epsilon-u^\epsilon+u^\epsilon\|_{H^1_{\epsilon}( \Omega^\epsilon)}^2 \\ \\ 
\qquad= \|P_{1+\eta}u^\epsilon-u^\epsilon\|_{H^1_{\epsilon}( \Omega^\epsilon)}^2 +
\|u^\epsilon\|_{H^1_{\epsilon}( \Omega^\epsilon)}^2 + 2 ( P_{1+\eta}u^\epsilon-u^\epsilon, u^\epsilon )_{H^1_{\epsilon}( \Omega^\epsilon)}\\ \\
\displaystyle\qquad= \|P_{1+\eta}u^\epsilon-u^\epsilon\|_{H^1_{\epsilon}( \Omega^\epsilon)}^2 +
\|u^\epsilon\|_{H^1_{\epsilon}( \Omega^\epsilon)}^2 + 2\int_{\Omega^\epsilon}(P_{1+\eta}u^\epsilon-u^\epsilon)f^\epsilon \, dx_1 dx_2 
\end{array}
$$
where we have used that  $u^\epsilon$ satisfies the variational formulation (\ref{VFP}) with  $\varphi=P_{1+\eta}u^\epsilon-u^\epsilon \in H^1(\Omega^\epsilon)$.

Consequently, it follows from (\ref{EQ1}) that
\begin{eqnarray*}
V_\epsilon(u^\epsilon) 
& \geq & \frac{1}{2}\|u^\epsilon\|_{H^1_\epsilon(\Omega^\epsilon\setminus \Omega^\epsilon(\frac{1}{1+\eta}))}^2+
\frac{1}{2(1+\eta)}\|P_{1+\eta}u^\epsilon-u^\epsilon\|_{H^1_{\epsilon}( \Omega^\epsilon)}^2 \\
& & +\frac{1}{2(1+\eta)} \| u^\epsilon \|_{H^1_{\epsilon}( \Omega^\epsilon)}^2  
- \frac{1}{1+\eta} \int_{\Omega^\epsilon}(u^\epsilon - P_{1+\eta}u^\epsilon)f^\epsilon \, dx_1 dx_2 
-\int_{\Omega^\epsilon} f^\epsilon u_\epsilon \, dx_1 dx_2 \\
& \geq & \frac{1}{2}\|u^\epsilon\|_{H^1_\epsilon(\Omega^\epsilon\setminus \Omega^\epsilon(\frac{1}{1+\eta}))}^2+
\frac{1}{2(1+\eta)}\|P_{1+\eta}u^\epsilon-u^\epsilon\|_{H^1_{\epsilon}( \Omega^\epsilon)}^2 \\
& & +\frac{1}{1+\eta}V_\epsilon(u^\epsilon) 
+ \int_{\Omega^\epsilon}\left(\frac{1}{1+\eta}P_{1+\eta}u^\epsilon-u^\epsilon \right) f^\epsilon \, dx_1 dx_2. 
\end{eqnarray*}

Hence, we obtain
\begin{equation} \label{EQ2}
\begin{gathered}
\frac{1}{2}\|u^\epsilon\|_{H^1_\epsilon(\Omega^\epsilon\setminus \Omega^\epsilon(\frac{1}{1+\eta}))}^2+
\frac{1}{2(1+\eta)}\|P_{1+\eta}u^\epsilon-u^\epsilon\|_{H^1_{\epsilon}( \Omega^\epsilon)}^2 \\
\qquad \quad \quad \leq \frac{\eta}{1+\eta}V(u^\epsilon) 
+ \int_{\Omega^\epsilon}\left(\frac{1}{1+\eta}P_{1+\eta}u^\epsilon-u^\epsilon \right) f^\epsilon \, dx_1 dx_2.
\end{gathered}
\end{equation}

Now, we analyze the integral 
$$
\int_{\Omega^\epsilon}\left(\frac{1}{1+\eta}P_{1+\eta}u^\epsilon-u^\epsilon \right) f^\epsilon \, dx_1 dx_2.
$$
To this, observe that
$$
u^\epsilon(x_1, x_2) - (P_{1+\eta} u^\epsilon)(x_1, x_2) = u^\epsilon(x_1, x_2) - u^\epsilon \left(x_1, \frac{x_2}{1+\eta} \right)
= \int_{\frac{x_2}{1+\eta}}^{x_2} \frac{\partial u^\epsilon}{\partial s}(x_1, s) ds
$$
which implies
\begin{eqnarray*}
\Big| u^\epsilon(x_1, x_2) - u^\epsilon \left(x_1, \frac{x_2}{1+\eta} \right) \Big| 
& \le & \int_{\frac{x_2}{1+\eta}}^{x_2} \Big| \frac{\partial u^\epsilon}{\partial s}(x_1, s) \Big| ds \\
& \le & \left( \int_{\frac{x_2}{1+\eta}}^{x_2} \Big| \frac{\partial u^\epsilon}{\partial s}(x_1, s) \Big|^2 ds \right)^{1/2} 
\left( \frac{\eta x_2}{1+\eta} \right)^{1/2}.
\end{eqnarray*}
Thus
$$
\begin{gathered}
\int_0^{G_\epsilon(x_1)} \Big| u^\epsilon(x_1, x_2) - u^\epsilon\left(x_1, \frac{x_2}{1+\eta}\right) \Big|^2 \, dx_2 \\
\le  \left( \int_0^{G_\epsilon(x_1)} \Big| \frac{\partial  u^\epsilon}{\partial s}(x_1, s) \Big|^2 ds \right) 
\left( \frac{\eta}{1+\eta} \right) \left( G_\epsilon(x_1) \right)^2  
 \end{gathered}
$$
and we have
$$
\| u^\epsilon - P_{1+\eta} u^\epsilon \|_{L^2(\Omega^\epsilon)} \le \left\| \frac{\partial  u^\epsilon}{\partial x_2} \right\|_{L^2(\Omega^\epsilon)}
\left( \frac{\eta}{1+\eta} \right)^{1/2} \, G_1.
$$
Consequently
\begin{eqnarray}
& & \left| \int_{\Omega^\epsilon} \left(\frac{1}{1+\eta}P_{1+\eta}u^\epsilon-u^\epsilon \right) f^\epsilon \, dx_1 dx_2 \right| \nonumber \\
& & \le \frac{\eta}{1+\eta} \int_{\Omega^\epsilon} \left| f^\epsilon u^\epsilon \right| \, dx_1 dx_2
+ \frac{1}{1+\eta} \int_{\Omega^\epsilon} \left| \left( P_{1+\eta}u^\epsilon-u^\epsilon \right) f^\epsilon \right|  \, dx_1 dx_2  \nonumber \\
& & \le \frac{\eta}{1+\eta} \| f^\epsilon \|_{L^2(\Omega^\epsilon)} \|u^\epsilon \|_{L^2(\Omega^\epsilon)}
+ \frac{\eta^{1/2}}{(1+\eta)^{3/2}} \| f^\epsilon \|_{L^2(\Omega^\epsilon)} 
\left\| \frac{\partial  u^\epsilon}{\partial x_2} \right\|_{L^2(\Omega^\epsilon)} G_1. \nonumber \\ \label{EQINT}
\end{eqnarray}
Then, it follows from (\ref{EQ2}) that
\begin{equation} \label{EQ3}
\begin{gathered}
\frac{1}{2}\|u^\epsilon\|_{H^1_\epsilon(\Omega^\epsilon\setminus \Omega^\epsilon(\frac{1}{1+\eta}))}^2 
+ \frac{1}{2(1+\eta)}\|P_{1+\eta}u^\epsilon-u^\epsilon\|_{H^1_{\epsilon}( \Omega^\epsilon)}^2  \\
 \leq \frac{\eta}{1+\eta} \left( V(u^\epsilon) + \| f^\epsilon \|_{L^2(\Omega^\epsilon)} \|u^\epsilon \|_{L^2(\Omega^\epsilon)} \right) 
 + \frac{\eta^{1/2}}{(1+\eta)^{3/2}} \| f^\epsilon \|_{L^2(\Omega^\epsilon)} 
\left\| \frac{\partial  u^\epsilon}{\partial x_2} \right\|_{L^2(\Omega^\epsilon)} G_1. 
\end{gathered}
\end{equation}
Hence, due to (\ref{EST0}), (\ref{ESTV}) and (\ref{EQ3}), we obtain
\begin{equation} \label{EQ4}
\begin{gathered}
\|u^\epsilon\|_{H^1_\epsilon(\Omega^\epsilon\setminus \Omega^\epsilon(\frac{1}{1+\eta}))}^2 
+ \|P_{1+\eta}u^\epsilon-u^\epsilon\|_{H^1_{\epsilon}( \Omega^\epsilon)}^2  \\
 \leq 2 \| f^\epsilon \|_{L^2(\Omega^\epsilon)}^2 
 \left( \frac{3}{2} \eta + \eta^{1/2} \, \epsilon \, G_1 \right). 
\end{gathered}
\end{equation}

On the other hand, we have 
\begin{eqnarray*}
V_\epsilon(u^\epsilon)
& = & \frac{1}{2}\|u^\epsilon\|_{H^1_\epsilon(\Omega^\epsilon)}^2 - \int_{\Omega^\epsilon}f^\epsilon u^\epsilon \, dx_1 dx_2  \\
& = & \frac{1}{2}\|P_{1+\eta} u^\epsilon\|_{H^1_{\epsilon,1+\eta}(\Omega^\epsilon({1+\eta}))}^2 
- \int_{\Omega^\epsilon}f^\epsilon u^\epsilon \, dx_1 dx_2  \\
& = & \frac{1}{2} \|P_{1+\eta} u^\epsilon\|_{H^1_{\epsilon,1+\eta}(\Omega^\epsilon(1+\eta) \setminus \Omega^\epsilon)}^2 +
\frac{1}{2}\|P_{1+\eta}u^\epsilon\|_{H^1_{\epsilon,1+\eta}( \Omega^\epsilon)}^2-\int_{\Omega^\epsilon}f^\epsilon u^\epsilon \, dx_1 dx_2 \\
& \geq & \frac{1}{2(1+\eta)} \left( \|P_{1+\eta} u^\epsilon\|_{H^1_{\epsilon}(\Omega^\epsilon(1+\eta) \setminus \Omega^\epsilon)}^2 +
\|P_{1+\eta}u^\epsilon\|_{H^1_{\epsilon}( \Omega^\epsilon)}^2 \right) - \int_{\Omega^\epsilon}f^\epsilon u^\epsilon \, dx_1 dx_2. 
\end{eqnarray*}
So, we can proceed as in (\ref{EQ4}) to get
\begin{equation} \label{EQ5}
\begin{gathered}
\|P_{1+\eta}u^\epsilon\|_{H^1_\epsilon(\Omega^\epsilon(1+\eta) \setminus \Omega^\epsilon)}^2 
+ \|P_{1+\eta}u^\epsilon-u^\epsilon\|_{H^1_{\epsilon}( \Omega^\epsilon)}^2  \\
 \leq 2 \| f^\epsilon \|_{L^2(\Omega^\epsilon)}^2 
 \left( \frac{3}{2} \eta + \eta^{1/2} \, \epsilon \, G_1 \right). 
\end{gathered}
\end{equation}
Putting together (\ref{EQ4}) and (\ref{EQ5}), we complete the proof of the lemma.
\end{proof}

\par\bigskip
We are in conditions now to proof the main result of this section.

\par\bigskip\noindent {\it Proof of Theorem \ref{BPT}.} 
If we define $\eta=\delta/G_0$, then under the hypotheses on $G_\eps$ and $\hat G_\eps$ and if $\|G_\eps-\hat G_\eps\|_{L^\infty(0,1)}\leq \delta$ we have that
\begin{equation}\label{SUB1}
\begin{array}{l}
\Omega^\eps (\frac{1}{1+\eta})\subset \hat\Omega^\eps\subset \Omega^\eps (1+\eta), \\
\hat\Omega^\eps (\frac{1}{1+\eta})\subset \Omega^\eps\subset \hat\Omega^\eps (1+\eta).
\end{array}
\end{equation}

Applying \eqref{SUB1} and  Lemma \ref{BSL} we easily get that 
\begin{equation}\label{first-estimate}
\begin{array}{l}
\|u^\eps\|^2_{H^1(\Omega^\eps\setminus\hat\Omega^\eps)}\leq 
\|u^\eps\|^2_{H^1(\Omega^\eps\setminus\Omega^\eps(\frac{1}{1+\eta}))}\leq C\sqrt{\eta}\leq C\sqrt{\delta}\\
\|\hat u^\eps\|^2_{H^1(\hat\Omega^\eps\setminus\Omega^\eps)}\leq 
\|\hat u^\eps\|^2_{H^1(\hat \Omega^\eps\setminus\hat \Omega^\eps(\frac{1}{1+\eta}))}\leq C\sqrt{\eta}\leq C\sqrt{\delta}.\\
\end{array}\end{equation}

So we will concentrate in the first term of \eqref{main-inequality}. 
Therefore, 
\begin{eqnarray*}
V_\epsilon(u^\epsilon) & \leq & V_\epsilon(\left(P_{1+\eta}\hat u^\epsilon\right)|_{\Omega^\epsilon}) \\
& = & \frac{1}{2} \|P_{1+\eta}\hat u^\epsilon \|_{H^1_\epsilon(\Omega^\epsilon)}^2 - \int_{\Omega^\epsilon} f^\epsilon (P_{1+\eta}\hat u^\epsilon )\, dx_1 dx_2 \\
& \leq & \frac{1}{2} \| P_{1+\eta}\hat u^\epsilon \|_{H^1_\epsilon (\hat \Omega^\epsilon(1+\eta))}^2 
- \int_{\hat \Omega^\epsilon} f^\epsilon P_{1+\eta}\hat u^\epsilon \, dx_1 dx_2 \\
& &\qquad\qquad+\int_{\hat \Omega^\epsilon \setminus \Omega^\epsilon} f^\epsilon P_{1+\eta}\hat u^\epsilon \, dx_1 dx_2. 
\end{eqnarray*}
But from (\ref{INPETA}) and (\ref{INPETA2}) we get
\begin{equation}\label{estimate1}
\|P_{1+\eta}\hat u^\epsilon \|_{H^1_\epsilon (\hat \Omega^\epsilon(1+\eta))}^2
\leq (1+\eta)\|\hat u^\epsilon\|_{H^1_\eps(\hat\Omega^\epsilon)}^2.
\end{equation}
Moreover
$$
\int_{\hat \Omega^\epsilon} f^\epsilon P_{1+\eta}\hat u^\epsilon \, dx_1 dx_2=
\int_{\hat \Omega^\epsilon} f^\epsilon (P_{1+\eta}\hat u^\epsilon -\hat u_\eps) \, dx_1 dx_2
+\int_{\hat \Omega^\epsilon} f^\epsilon \hat u^\epsilon \, dx_1 dx_2$$
but from Lemma \ref{BSL},
$$
|\int_{\hat \Omega^\epsilon} f^\epsilon (P_{1+\eta}\hat u^\epsilon -\hat u_\eps) \, dx_1 dx_2|
\leq \|f^\epsilon\|_{L^2(\hat\Omega^\eps)}\|P_{1+\eta}\hat u^\epsilon -\hat u_\eps\|_{L^2(\hat\Omega^\eps)}
\leq C{\eta}^{1/4}.
$$
Also, 
$$|\int_{\hat \Omega^\epsilon \setminus \Omega^\epsilon} f^\epsilon P_{1+\eta}\hat u^\epsilon \, dx_1 dx_2|
\leq \|f^\epsilon\|_{L^2(\hat\Omega_\eps)}\|P_{1+\eta}\hat u^\eps\|_{L^2(\hat\Omega_\eps\setminus\Omega_\eps)}\leq C\eta^{1/4}$$
where we have used Lemma \ref{BSL} and \eqref{SUB1}.

Hence, putting all this information together, we get
$$
V_\epsilon(u^\epsilon)\leq (1+\eta)\hat V_\eps (\hat u^\eps)+\eta\|f^\eps\|_{L^2(\hat\Omega_\eps)}\|\hat u^\eps\|_{L^2(\hat\Omega^\eps)}+C\eta^{1/4}\leq \hat V_\eps(\hat u_\eps)+C\eta^{1/4}.
$$

Therefore
\begin{equation} \label{EQ6}
V_\epsilon(u^\epsilon) \leq \hat V_\epsilon(\hat u^\epsilon) + C \delta^{1/4}.
\end{equation}

On the other hand, we obtain by (\ref{NPETA}), (\ref{INPETA}) and \eqref{SUB1}
\begin{eqnarray}
V_\epsilon(u^\epsilon) & = & 
\frac{1}{2} \| u^\epsilon \|_{H^1_\epsilon(\Omega^\epsilon)}^2 - \int_{\Omega^\epsilon} f^\epsilon  u^\epsilon \, dx_1 dx_2  \nonumber \\ 
& = & \frac{1}{2}\|P_{1+\eta} u^\epsilon \|^2_{H^1_{\epsilon,1+\eta}(\Omega^\epsilon(1+\eta))} - \int_{\Omega^\epsilon} f^\epsilon u^\epsilon \, dx_1 dx_2  \nonumber \\
& \geq & \frac{1}{2(1+\eta)} \| P_{1+\eta}u^\epsilon \|^2_{H^1_\epsilon(\hat \Omega^\epsilon)} - \int_{\Omega^\epsilon} f^\epsilon u^\epsilon \, dx_1 dx_2  \nonumber \\
& = & \frac{1}{2(1+\eta)} \|P_{1+\eta}u^\epsilon -\hat u^\epsilon + \hat u^\epsilon \|^2_{H^1_\epsilon(\hat \Omega^\epsilon)} - \int_{\Omega^\epsilon}f^\epsilon u^\epsilon \, dx_1 dx_2  \nonumber \\
& = & \frac{1}{2(1+\eta)} \left(\|P_{1+\eta}u^\epsilon - \hat u^\epsilon \|^2_{H^1_\epsilon(\hat \Omega^\epsilon)} + 
\|\hat u^\epsilon\|^2_{H^1_\epsilon(\hat \Omega^\epsilon)}\right.  \nonumber \\
& & \left. +
2(P_{1+\eta}u^\epsilon - \hat u^\epsilon,\hat u^\epsilon)_{H^1_\epsilon(\hat \Omega^\epsilon)}\right) 
- \int_{\Omega^\epsilon} f^\epsilon u^\epsilon \, dx_1 dx_2  \nonumber \\
& = & \frac{1}{2(1+\eta)} \left(\|P_{1+\eta}u^\epsilon - \hat u^\epsilon \|^2_{H^1_\epsilon(\hat \Omega^\epsilon)} + 
\|\hat u^\epsilon\|^2_{H^1_\epsilon(\hat \Omega^\epsilon)} \right. \nonumber \\ 
& & \left. + 2 \int_{\hat \Omega^\epsilon} (P_{1+\eta}u^\epsilon - \hat u^\epsilon) f^\epsilon \, dx_1 dx_2 \right) 
- \int_{\Omega^\epsilon} f^\epsilon u^\epsilon \, dx_1 dx_2  \nonumber \\
& = & \frac{1}{2(1+\eta)}\|P_{1+\eta}u^\epsilon - \hat u^\epsilon \|^2_{H^1_\epsilon(\hat \Omega^\epsilon)} 
+ \frac{1}{1+\eta} \hat V_\epsilon(\hat u^\epsilon) \nonumber \\
& & + \frac{1}{1+\eta} \int_{\hat \Omega^\epsilon} P_{1+\eta}u^\epsilon \, f^\epsilon  \, dx_1 dx_2 
- \int_{\Omega^\epsilon} f^\epsilon u^\epsilon  \, dx_1 dx_2. \label{EQ7}
\end{eqnarray}

But, due to (\ref{EQINT}) and Lemma \ref{BSL}, we have
\begin{eqnarray*}
& & \Big| \frac{1}{1+\eta} \int_{\hat \Omega^\epsilon} P_{1+\eta}u^\epsilon \, f^\epsilon  \, dx_1 dx_2 
- \int_{\Omega^\epsilon} f^\epsilon u^\epsilon  \, dx_1 dx_2 \Big| \\
& \le & \frac{1}{1+\eta} \Big| \int_{\hat \Omega^\epsilon \setminus \Omega^\epsilon} P_{1+\eta}u^\epsilon \, f^\epsilon  \, dx_1 dx_2  \Big| \\
& & + \Big| \frac{1}{1+\eta} \int_{ \Omega^\epsilon} P_{1+\eta}u^\epsilon \, f^\epsilon  \, dx_1 dx_2 
- \int_{\Omega^\epsilon} f^\epsilon u^\epsilon  \, dx_1 dx_2  \Big|  \\
& & \le \frac{1}{1+\eta} \left( \| P_{1+\eta}u^\epsilon \|_{L^2(\hat \Omega^\epsilon \setminus \Omega^\epsilon)} 
\| f^\epsilon \|_{L^2(\hat \Omega^\epsilon \setminus \Omega^\epsilon)} +  C_0 \, \eta^{1/2} \right) \\
& & \le \frac{1}{1+\eta} \left( \| P_{1+\eta}u^\epsilon \|_{L^2(\Omega^\epsilon(1+\eta) \setminus \Omega^\epsilon)} 
\| f^\epsilon \|_{L^2( \Omega^\epsilon(1+\eta) \setminus \Omega^\epsilon)} +  C_0 \, \eta^{1/2} \right) \\
& & \le \frac{1}{1+\eta} \left( \| P_{1+\eta}u^\epsilon \|_{H^1_\epsilon(\Omega^\epsilon(1+\eta) \setminus \Omega^\epsilon)} 
\| f^\epsilon \|_{L^2( \Omega^\epsilon(1+\eta) \setminus \Omega^\epsilon)} +  C_0 \, \eta^{1/2} \right) \\
& & \le \frac{1}{1+\eta} \left( C_1 \, \eta^{1/4} \| f^\epsilon \|_{L^2( \Omega^\epsilon(1+\eta) \setminus \Omega^\epsilon)} 
+  C_0 \, \eta^{1/2} \right) \le C_2 \, \delta^{1/4} 
\end{eqnarray*}
with $C_0$, $C_1$ and $C_2$ independent of $\epsilon$ and $\eta$.

Consequently
\begin{equation} \label{EQ8}
V_\epsilon(u^\epsilon) \geq  
\frac{1}{1+\eta} \hat V_\epsilon (\hat u^\epsilon) + \|P_{1+\eta}u^\epsilon 
- \hat u^\epsilon \|^2_{H^1_\epsilon(\hat \Omega^\epsilon)} - C_ 2 \, \delta^{1/4}.
\end{equation}
Thus, it follows from (\ref{EQ6}) and (\ref{EQ8}) that
$$
\|P_{1+\eta}u^\epsilon - \hat u^\epsilon \|^2_{H^1_\epsilon(\hat \Omega^\epsilon)} 
\leq \frac{\eta}{1+\eta} \hat V_\epsilon(\hat u^\epsilon) + C \, \delta^{1/4}  + C_ 2 \, \delta^{1/4}
$$
which implies 
\begin{equation}\label{almost-last-one}
\|P_{1+\eta}u^\epsilon - \hat u^\epsilon \|^2_{H^1_\epsilon(\hat \Omega^\epsilon)} 
\leq C_3 \, \delta^{1/4}
\end{equation}
with $C_3$ independent of $\epsilon$ and $\eta$.

But, from Lemma \ref{BSL} we have that $\|u^\eps-P_{1+\eta}u^\eps\|_{H^1_\eps(\Omega^\eps)}^2\leq C\sqrt{\eta}$. 
Hence, from \eqref{almost-last-one} we get 
\begin{equation} \label{EQ001}
\|u^\epsilon - \hat u^\epsilon \|^2_{H^1_\epsilon( \Omega^\epsilon\cap\hat\Omega_\eps)}  \leq C \, \delta^{1/4}
\end{equation}
for $C$ independent of $\epsilon$. This provers the result.
\cqd

\section{The General Case}
\label{general-case}
We provide now a proof of the main result, Theorem \ref{GT}.

\begin{proof}
From estimate (\ref{EST00})  and (\ref{EQOP}) we derive 
$$
\begin{gathered}
\| P_{\epsilon} u^\epsilon \|_{L^2(\hat \Omega)}, \Big\| \frac{\partial P_{\epsilon} u^\epsilon}{\partial x_1} \Big\|_{L^2(\hat\Omega)} \textrm{ and }
\frac{1}{\epsilon} \Big\| \frac{\partial P_{\epsilon} u^\epsilon}{\partial x_2} \Big\|_{L^2(\hat\Omega)} \le M 
\textrm{ for all } \epsilon > 0
\end{gathered}
$$
where $M$ is a positive constant independent of $\epsilon$.
Hence, there exists $u_0\in H^1(\hat \Omega)$  and a subsequence, still denoted by $P_{\epsilon} u^\epsilon$, such that
\begin{equation} \label{LEO-bis}
\begin{gathered}
P_{\epsilon} u^{\epsilon} \rightharpoonup u_0 \quad w-H^1(\hat \Omega) \\
P_{\epsilon} u^{\epsilon} \rightarrow u_0 \quad s-L^2(\hat\Omega) \\
\textrm{ and }
\frac{\partial P_{\epsilon} u^\epsilon}{\partial x_2} \rightarrow 0 \quad s-L^2(\hat \Omega).
\end{gathered}
\end{equation}

It follows from (\ref{LEO-bis}) that $u_0(x_1,x_2) = u_0(x_1)$ on $\hat \Omega$. 

We will show that $u_0$ satisfies the Neumann problem (\ref{VFPDL}).
To do this, we use a discretization argument.

Let us fix a parameter $\delta>0$ small enough and consider a function $G^\delta (x,y)$ with the property that 
$0\leq G^\delta(x,y)-G(x,y)\leq \delta$ for all $(x,y)\in I\times \R$ and such that the function $G^\delta$ satisfies hypothesis
(H) and is piecewise periodic
as described in Section \ref{piecewise-periodic}.  To construct this function, we may proceed as follows. The function $G$ is uniformly $C^1$ in each of the domains $(\xi_{i-1},\xi_i)\times \R$ and it is also periodic in the second variable. In particular, for $\delta>0$ small enough and for a fixed $z\in (\xi_{i-1},\xi_i)$ we have that there exists a small interval $(z-\eta,z+\eta)$ with $\eta$ depending only on $\delta$ such that  $|G(x,y)-G(z,y)|+|\partial_y G(x,y)-
\partial_y G(z,y)|<\delta/2$ for all $x\in (z-\eta,z+\eta)\cap (\xi_{i-1},\xi_i)$ and for all $y\in \R$. This allows us to select a finite number of points:  
$\xi_{i-1}=\xi_{i-1}^1< \xi_{i-1}^2<\ldots<\xi_{i-1}^r=\xi_i$ such that $\xi_{i-1}^r-\xi_{i-1}^{r-1}<\eta$ and therefore, defining
$G^\delta(x,y)=G(\xi_{i-1}^r,y)+\delta/2$ for all $x\in (\xi_{i-1}^r,\xi_{i-1}^{r+1})$ we have that $0\leq G^\delta(x,y)-G(x,y)\leq \delta$ and $|\partial_y G^\delta(x,y)-\partial_y G(x,y)|\leq \delta$ for all
$(x,y)\in (\xi_{i-1},\xi_i)\times \R$. This construction can be done for all $i=1,\ldots, N$.

 In particular, if we rename all the points $\xi_i^k$ constructed above by $0=z_0<z_1<\ldots<z_m=1$ (and observe that $m=m(\delta)$),  the function $G^\delta(x,y)=G^\delta_i(y)$ with $(x,y)\in (z_{i-1},z_i)\times \R$, $i=1,\ldots, m$ and $G^\delta_i$ is $C^1$ and $L$-periodic.

We denote by $G_\eps^\delta(x)=G^\delta(x,x/\eps)$ and consider the domains
$$
\begin{gathered}
\Omega^{\epsilon,\delta} = \{ (x,y) \in \R^2 \; | \;  x \in I,  \quad
 0 < y < G_{\epsilon}^{\delta}(x) \} \\
\Omega^\delta = I \times (0, G_1)\\ 
\hat \Omega^\delta = \Omega^\delta \setminus 
\cup_{i=1}^{m-1} \{ (z_i, y) \in \R^2 \, | \,  G_0  < y <G_1\}.
\end{gathered}
$$

Since $G_{\epsilon}^{\delta}$ satisfies the hyphoteses of the Lemma \ref{EOT}, there exists an extension operator
$$
P_{\epsilon,\delta} \in \mathcal{L}(L^p(\Omega^{\epsilon,\delta}),L^p(\hat\Omega^\delta)) 
  \cap \mathcal{L}(W^{1,p}(\Omega^{\epsilon,\delta}),W^{1,p}(\hat\Omega^\delta))
$$
satisfying the uniform estimate (\ref{EQOP}) with $\eta(\eps)\sim 1/\eps$.

Since $f^\eps\in L^2(\Omega_\eps)$ with $\|f^\eps\|_{L^2(\Omega^\eps)}\leq C$, then if we extend this function
by 0 outside $\Omega^\eps$ and denoting the extended function again by $f^\eps$ and using that $G_\delta\geq G$, we have that 
$\hat f^\eps_\delta(x)=\int_0^{G_\eps^\delta(x)}f^\eps (x,y)dy=\int_0^{G_\eps(x)}f^\eps(x,y)dy=\hat f^\eps(x)$ and
by hypothesis, we have that $\hat f^\eps_\delta\equiv \hat f^\eps\rightharpoonup \hat f$ w-$L^2(0,1)$. 

Therefore, it follows from Theorem \ref{PPCT} that for each $\delta > 0$ fixed, 
there exist $u^\delta \in H^1(0,1)$
such that the solutions $u^{\epsilon,\delta}$ 
of (\ref{P}) in $\Omega^{\epsilon,\delta}$ satisfy
\begin{equation} \label{UDC}
P_{\epsilon,\delta} u^{\epsilon,\delta} \rightharpoonup u^\delta \quad w - H^1(\hat \Omega^\delta)
\end{equation}
where $u^\delta(x_1,x_2) = u^\delta(x_1)$ on $\Omega^\delta$ is the unique solution of the Neumann problem

\begin{equation} \label{VFPD}
\int_{I} \Big\{ r^\delta(x) \; u_x^\delta(x) \, \varphi_x(x) 
+ p^\delta(x) \, u^\delta(x) \, \varphi(x) \Big\} dx = \int_{I}  \, \hat f(x) \, \varphi(x) \, dx, \quad \forall \varphi \in H^1(I),
\end{equation}
where $r^\delta: I \mapsto \R$ and $p^\delta: I \mapsto \R$ are strictly positive functions, 
locally constant, given by
\begin{equation} \label{RPFD}
\left\{ 
\begin{gathered}
r^\delta(x) = r_{i,\delta}(x) \\
p^\delta(x) = \frac{|Y^*_i|}{L} 
\end{gathered}
\right. 
\quad x \in (z_{i-1},z_i)
\end{equation}
where 
\begin{equation}\label{def-rdelta}
r_{i,\delta} = \frac{1}{L} \int_{Y^*_i} \Big\{ 1 - \frac{\partial X_i}{\partial y_1}(y_1,y_2) \Big\} dy_1 dy_2.\
\end{equation}
The function $X_i$ is the unique solution of (\ref{AUX}) in the representative cells $Y^*_i$ defined in (\ref{CELL}) 
for all $i = 1, ..., m$.

Now, we pass to the limit in $(\ref{VFPD})$ as $\delta \to 0$.
To do this, we consider the functions $r^\delta$ and $p^\delta$ defined in $x \in I$
and the functions $r$ and $p$ defined in (\ref{RPFL}).
We have that $r^\delta$ and $p^\delta$ converge to $r$ and $p$ uniformly in $I$.
The uniform convergence of $r^\delta$ to $r$ in $I$ it follows from  Proposition \ref{prop-appendix} proved in the Appendix.
The uniform convergence of $p^\delta$ to $p$ follows from the uniform convergence of $G^\delta$ to $G$ as
$\delta\to 0$. 

Since we have the uniform convergence of $r^\delta$ and $p^\delta$ to $r$ and $p$ respectively,
we have by \cite[p. 8]{BLP} or \cite[p. 1]{CP} the following limit variational formulation:
to find $u \in H^1(I)$ such that 
\begin{equation} \label{VFPDL-F}
\int_{I} \Big\{ r(x) \; u_x(x) \, \varphi_x(x) 
+ p(x) \, u(x) \, \varphi(x) \Big\} dx = \int_{I}  \, \hat f(x) \, \varphi \, dx
\end{equation}
for all $\varphi \in H^1(I)$.

Hence, there exists $u^* \in H^1(I)$ such that 
\begin{equation} \label{DCU}
u^\delta \to u^* \textrm{ in } H^1(I)
\end{equation}
where $u^*$ is the unique solution of the Neumann problem (\ref{VFPDL-F}).

To finish the proof, we need to show that $u^* = u_0$ in $I$, where $u_0$ is the function obtained in (\ref{LEO-bis}).

Let us consider the open square $\Omega_0=I\times (0,G_0)$ which satisfies $\Omega_0\subset \Omega_\delta^\eps$ for all
$\delta$ and $\eps$ small enough. Observe that $\|u^*-u_0\|_{L^2(I)}^2=G_0^{-1}\|u^*-u_0\|_{L^2(\Omega_0)}^2$ and therefore, to show that $u^*=u_0$ it is enough to show that $\|u^*-u_0\|_{L^2(\Omega_0)}^2=0$. Hence, adding and
substracting the appropriate functions and with the triangular inequality, 
\begin{equation}\label{EQF001}
\begin{array}{l}
\| u^* - u_0 \|_{L^2(\Omega_0)} \leq 
\| u^* - u^{\delta}  \|_{L^2(\Omega_0)}+ 
\| u^\delta - u^{\epsilon,\delta}  \|_{L^2(\Omega_0)}\\ \\
\qquad \qquad+ \|  u^{\epsilon,\delta} - u^\epsilon   \|_{L^2(\Omega_0)} 
+ \|  u^\epsilon - u_0  \|_{L^2(\Omega_0)} \end{array}
\end{equation}
for all $\epsilon$ and $\delta > 0$.

Now, let $\eta$ be a positive small number. 
From \eqref {DCU} and from Theorem \ref{BPT}, we can choose a $\delta>0$ fixed and small such that 
$\| u^* - u^{\delta}  \|_{L^2(\Omega_0)}\leq \eta$ and $\|  u^{\epsilon,\delta} - u^\epsilon   \|_{L^2(\Omega_0)} \leq \eta$
uniformly for all $\eps>0$.  Now, from \eqref{UDC} for this particular value of $\delta$, we can choose $\epsilon_1>0$ small
enough such that $\| u^\delta - u^{\epsilon,\delta}  \|_{L^2(\Omega_0)}\leq \eta$  for $0<\eps<\eps_1$. Moreover, from 
\eqref{LEO-bis} we have that there exists $\eps_2>0$ such that $ \|  u^\epsilon - u_0  \|_{L^2(\Omega_0)}\leq \eta$ for
all $0<\eps<\eps_2$.   Hence with $\eps=\min\{\eps_1,\eps_2\}$ applied to \eqref{EQF001}, we get
$\| u^* - u_0 \|_{L^2(\Omega_0)}\leq 4\eta$. Since $\eta$ is arbitrarilly small, then $u^*=u_0$. 
\end{proof}

\appendix
\section{ A perturbation result in the basic cell.}
\label{appendix}

In the proof of the main result in Section \ref{general-case} we have used the convergence of $r^\delta\to r$, where
$r^\delta=r_{i,\delta}$ in the interval $(z_{i-1},z_i)$, and $r_{i,\delta}$, $r$ are defined by (\ref{def-rdelta}) and (\ref{RPFL}), respectively.  In order to proof such a convergence we need to analyze how the function $X$, solution of
\begin{equation} \label{AUXKato1}
\left\{
\begin{gathered}
- \Delta X = 0  \textrm{ in } Y^* \\
\frac{\partial X}{\partial N} = 0  \textrm{ on } B_2 \\
\frac{\partial X}{\partial N} = N_1 = - \frac{\partial_2 G}{\sqrt{1 + \left( \partial_2 G \right)^2}}  \textrm{ on } B_1 \\
X(0,y_2) = X(L,y_2) \textrm{ on } B_0 \\
\int_{Y^*} X \; dy_1 dy_2 = 0
\end{gathered}
\right.
\end{equation}
on the representative cell
\begin{equation} \label{DKATO}
Y^* = \{ (y_1,y_2) \in \R^2 \; | \; 
0 < y_1 < L, \; \; 0 < y_2 < G(y_1) 
\},
\end{equation}
depends on the function $G$.

We will consider the following class of admissible functions
\begin{equation}\label{admissible}
A(M)=\{G\in C^1(\R), L-\hbox{periodic},\, 0<G_0\leq G(\cdot)\leq G_1, |G'(s)|\leq M\}
\end{equation}
and we will denote by $Y^*(G)$ and $X(G)$ the basic cell (\ref{DKATO}) and the 
 solution of (\ref{AUXKato1}) for a particular  $G\in A(M)$. 
Observe that for each $G\in A(M)$, we have the extension operator $E_G:H^1(Y^*(G))\to H^1(Y)$ 
as constructed in Lemma \ref{EOT}, which will satisfy $\|E_G u\|_{H^1(Y)}\leq C\|u\|_{H^1(Y^*(G))}$
with $C=C(M)$, but $C$ independent of $G$, and therefore this constant  can be chosen the same for all $G\in A(M)$.

For each $\hat G \in A(M)$, we can consider the basic cell $Y^*(\hat G)$ defined by (\ref{DKATO}). 
Here, we proceed as \cite[p. 84]{JA} and \cite{HR}, we begin by making the following transformation on the domain $Y^*(G)$
\begin{eqnarray*}
L_{\hat G} & : & Y^*(G) \mapsto Y^*(\hat G) \\
& & (z_1, z_2) \to (z_1, \hat F(z_1) \; z_2) = (y_1, y_2)
\end{eqnarray*}
where $\hat F(z) = \frac{\hat G(z)}{G(z)}$.
The Jacobian matrix for this transformation is
$$
J L_{\hat G} (z_1,z_2) = \left(
\begin{array}{cc}
1 & 0 \\
\hat F'(z_1) z_2 & \hat F(z_1)
\end{array}
\right)
$$
and observe that its determinant is given by $|J L_{\hat G} (z_1,z_2)|=\hat F(z_1)=\frac{\hat G(z)}{G(z)}$. 

Using $L_{\hat G}$, we can show that problem (\ref{AUXKato1}) in $Y^*(\hat G)$ is equivalent to
\begin{equation} \label{AUXKato00}
\left\{
\begin{array}{ll}
- \frac{1}{\hat F} {\rm div }( B_{\hat G} W ) = 0 & \textrm{ in } Y^*(G) \\
B_{\hat G} W \cdot N = 0 & \textrm{ on } B_2(G) \\
B_{\hat G} W \cdot N = - \frac{\partial_2 \hat G}{\sqrt{1 + \left( \partial_2 \hat G \right)^2}} & \textrm{ on } B_1(G) \\
W(0,z_2) = W(L,z_2) & \textrm{ on } B_0(G) \\
\int_{Y^*(G)} W \, |J L_{\hat G}| \, dz_1 dz_2 = 0
\end{array}
\right. 
\end{equation}
where 
$$
(B_{\hat G} W)(z_1,z_2) = \left(
\begin{array}{c}
\hat F(z_1) \; \frac{\partial U}{\partial z_1}(z_1,z_2) - \hat F'(z_1) \; z_2 \; \frac{\partial U}{\partial z_2}(z_1,z_2) \\
- \hat F'(z_1) \; z_2 \; \frac{\partial U}{\partial z_1}(z_1,z_2) 
+ \frac{1}{\hat F(z_1)} \; ( 1 + (z_2 \; \hat F'(z_1))^2 ) \; \frac{\partial U}{\partial z_2}(z_1,z_2)
\end{array}
\right).
$$
That is, we have that $X(\hat G)$ is the solution of (\ref{AUXKato1}) in $Y^*(\hat G)$ 
if and only if $W(\hat G) = X(\hat G) \circ L_{\hat G}$ satisfies equation (\ref{AUXKato00}) in $Y^*(G)$.

Moreover, if we define $$U(\hat G)=W(\hat G)-\frac{1}{|Y^*(G)|}\int_{Y^*(G)}W(\hat G)$$
then, $U(\hat G)$ is the unique solution of 
\begin{equation} \label{AUXKato0}
\left\{
\begin{array}{ll}
- \frac{1}{\hat F} {\rm div }( B_{\hat G} U) = 0 & \textrm{ in } Y^*(G) \\
B_{\hat G} U \cdot N = 0 & \textrm{ on } B_2(G) \\
B_{\hat G} U \cdot N = - \frac{\partial_2 \hat G}{\sqrt{1 + \left( \partial_2 \hat G \right)^2}} & \textrm{ on } B_1(G) \\
U(0,z_2) = U(L,z_2) & \textrm{ on } B_0(G) \\
\int_{Y^*(G)} U  \, dz_1 dz_2 = 0
\end{array}
\right. .
\end{equation}

Also, we have that the variational formulation of  problem (\ref{AUXKato0}) 
is given by the bilinear form 
\begin{equation} \label{b0}
\begin{array}{rl}
\rho_{\hat G}: H(G) \times H(G) &\mapsto \R \\
 \displaystyle (U,V) &\mapsto \displaystyle  \int_{Y^*(G)} B_{\hat G} U \cdot \nabla \left( \frac{V}{\hat F} \right) \; dz_1 dz_2 \nonumber 
\end{array}
\end{equation}
where $H(G)$ is given by 
\begin{eqnarray*}
H(G) & = & \{ U \in H^1(Y^*(G)) \, | \, U(0,z_2) = U(L,z_2) \textrm{ on } B_0(G), \, \int_{Y^*(G)}U=0 \}
\end{eqnarray*}
with the $H^1(Y^*(G))$ norm. 
Hence, we have that $X(\hat G)$ is the solution of (\ref{AUXKato1}) in $Y^*(\hat G)$ if and only if 
\begin{equation}\label{def-of-U}
U(\hat G) = X(\hat G) \circ L_{\hat G}-\frac{1}{|Y^*(G)|}\int_{Y^*(G)}X(\hat G) \circ L_{\hat G}
\end{equation}
satisfies 
\begin{equation}\label{weak-form-hat G}
\rho_{\hat G}(U,V) = 
- \int_{B_1(G)} \frac{\hat G'}{\sqrt{1 + ( \hat G ')^2}} \; \frac{V}{\hat F} \; d\mathcal{S},  \quad
\textrm{ for all } V \in H(G).
\end{equation}

Observe also that the bilinear form associated to problem (\ref{AUXKato1}) in $Y^*(G)$ is given
by 
\begin{equation} \label{b0G}
\begin{array}{rl}
\rho_{G}: H(G) \times H(G) &\mapsto \R \\
 \displaystyle (U,V) &\mapsto \displaystyle  \int_{Y^*(G)} \nabla  U \cdot \nabla V  \; dz_1 dz_2 \nonumber 
\end{array}
\end{equation}
and the weak formulation is given by 
$$
\rho_{G}(U,V) = 
- \int_{B_1(G)} \frac{G'}{\sqrt{1 + ( G ')^2}} \; V \; d\mathcal{S}, \quad 
\textrm{ for all } V \in H(G).
$$

As a matter of fact, we will be able to show the following: 

\begin{proposition}\label{prop-appendix}
Let us consider the family of admissible functions $G \in A(M)$ for some constant $M$, where $A(M)$ is 
defined in (\ref{admissible}). 

Then, for each $\eps>0$ there exists $\delta>0$  such that if $G, \hat G \in A(M)$ with 
$\|G - \hat G \|_{C^1(\R)}\leq \delta$, then
\begin{equation}\label{conv-h1}
\| U(G) - U(\hat G) \|_{H(Y^*(G))} \leq \epsilon.
\end{equation}

Moreover, we obtain from (\ref{conv-h1})
\begin{equation}\label{conv-r}
|r(G)-r(\hat G)|\leq \eps
\end{equation}
where 
\begin{equation}\label{def-rG}
r(G)=\frac{1}{|Y^*(G)|}\int_{Y^*(G)}
\Big\{ 1 - \frac{\partial X(G)}{\partial y_1}(y_1,y_2) \Big\} dy_1 dy_2
\end{equation}
and similarly for $r(\hat G)$. 
\end{proposition}

\begin{proof}
Obviously (\ref{conv-r}) follows straightforward from (\ref{conv-h1}), the relation between $U$ and $X$ stated
in  (\ref{def-of-U}), the definition of $r(G)$ in (\ref{def-rG})
and the fact that $G$ and $\hat G$ are close in the $C^1$-metric, and in particular in the uniform norm. 

Hence, we need to show (\ref{conv-h1}). For this
let $\rho_{G}$ and $\rho_{\hat G}$ be bilinear forms associated to the variational formulation of (\ref{AUXKato0}) for $G$ and $\hat G$ respectively.
We can get the following estimate
\begin{eqnarray}
& & | \rho_{\hat G}(U,V) - \rho_{G}(U,V) |   \nonumber  \\
& \le & \int_{Y^*(G)} \Big\{ 
\Big| \frac{\partial}{\partial z_1} \left( \frac{V}{\hat F} \right) \Big( (\hat F - 1) \frac{\partial U}{\partial z_1} 
 - z_2 \hat F' \frac{\partial U}{\partial z_2} \Big) \Big|  \nonumber \\
& & + \Big| \frac{1}{\hat F} \, \frac{\partial V}{\partial z_2} \Big( 
\frac{1}{\hat F} \Big( 1- \hat F + ( z_2 \hat F')^2 \Big)\frac{\partial U}{\partial z_2} 
- z_2 \hat F' \frac{\partial U}{\partial z_1}
\Big) \Big| \Big\} dz_1 dz_2  \nonumber \\
& \le & \frac{\| \hat F - 1 \|_{L^\infty}}{G_0} 
\Big\| \frac{\partial V}{\partial z_1} \Big\|_{L^2} \Big\| \frac{\partial U}{\partial z_1} \Big\|_{L^2} 
 \nonumber 
 + \frac{\| 1 - \hat F \|_{L^\infty} + G_1^2 \, \| \hat F' \|_{L^\infty}}{G_0^2} 
 \Big\| \frac{\partial V}{\partial z_2} \Big\|_{L^2} \Big\| \frac{\partial U}{\partial z_2}  \Big\|_{L^2} \nonumber \\
& &  + \frac{G_1}{G_0} \, \| \hat F' \|_{L^\infty} \Big( \Big\| \frac{\partial V}{\partial z_1} \Big\|_{L^2} 
\Big\| \frac{\partial U}{\partial z_2}  \Big\|_{L^2} + \Big\| \frac{\partial V}{\partial z_2} \Big\|_{L^2} \Big\| \frac{\partial U}{\partial z_1}  \Big\|_{L^2} \Big) \nonumber \\
& & + \frac{\| \hat F' \|_{L^\infty}}{G_0^2} \left( \| 1 - \hat F \|_{L^\infty} \Big\| \frac{\partial U}{\partial z_1}  \Big\|_{L^2}
+ G_1 \, \| \hat F' \|_{L^\infty} \Big\| \frac{\partial U}{\partial z_2}  \Big\|_{L^2} \right) \Big\| V \Big\|_{L^2}. 
\label{EKato} 
\end{eqnarray}
Since
$$
\begin{gathered}
\hat F(z_1) - 1 = \frac{\hat G(z_1) - G(z_1)}{G(z_1)}  \textrm{ and } \\
\hat F'(z_1) = \frac{ G'(z_1) \left( G(z_1) - \hat G(z_1) \right) 
+ G(z_1) \left( \hat G'(z_1) - G'(z_1) \right)}{G(z_1)^2} 
\end{gathered}
$$ 
we have
$$
\begin{gathered}
\| \hat F - 1 \|_{L^\infty( (0,L))} \le  
\frac{1}{G_0} \| G - \hat G  \|_{L^\infty( (0,L))} \textrm{ and } \\
\| \hat F' \|_{L^\infty((0,L))} \le  \frac{1}{{G_0}^2} 
\left( M \, \| G - \hat G \|_{L^\infty((0,L))} + G_1 \, \| G' - \hat G' \|_{L^\infty(0,L)} \right).
\end{gathered}
$$
Then, due to (\ref{EKato}), we get 
\begin{equation} \label{EK1}
| \rho_{\hat G}(U,V) - \rho_{G}(U,V) | \le  
 C  \left( \| G - \hat G \|_{L^\infty} + \| G' - \hat G' \|_{L^\infty} \right) 
\| U \|_{H} \| V \|_{H},
\end{equation}
where the constant $C$ depends on $M$, $G_0$ and $G_1$ and therefore it can be chosen the same
constant for all $G,\hat G\in A(M)$. 

Now, since $\rho_G$ is a coercive bilinear form in $H(G)$ (observe that we have imposed the condition 
$\int_{Y^*(G)}U=0$ in $H(G)$) then, there will exist a constant $c>0$ such that  $\rho_G(U(G)- U\hat (G),U(G)-U(\hat G))\geq c\|U(G)-U(\hat G)\|_{H(G)}^2$.  But, to simplify, if we denote by $U=U(G)$ and $\hat U=U(\hat G)$, we will have
$$ 
\begin{array}{l}
c\|U-\hat U\|_{H(G)}^2\leq \rho_G(U-\hat U, U-\hat U)= \rho_G(U, U-\hat U)-
\rho_G(\hat U, U-\hat U)\\ \\
\qquad \leq |\rho_G(U, U-\hat U)-\rho_{\hat G}(\hat U, U-\hat U)|+|\rho_{\hat G}(\hat U, U-\hat U)-\rho_G(\hat U, U-\hat U) |\\ \\
\displaystyle \leq \int_{B_1(G)}\left| \frac{G'}{\sqrt{1 + ( G ')^2}} - \frac{\hat G'}{\sqrt{1 + ( \hat G ')^2}} \; \frac{1}{\hat F}\right| |U-\hat U| \; d\mathcal{S}+C\|G-\hat G\|_{C^1(\R)}\|\hat U\|_{H(G)} \|U-\hat U|\|_{H(G)} .
\end{array}
$$

But using appropriate trace theorems in $B_1(G)$, the fact that $G$, $\hat G\in A(M)$ so that  they are uniformly bounded in $C^1(\R)$
(and therefore $\left| \frac{G'}{\sqrt{1 + ( G ')^2}} - \frac{\hat G'}{\sqrt{1 + ( \hat G ')^2}} \; \frac{1}{\hat F}\right|\leq C\|G-\hat G\|_{C^1(\R)}$ for a constant $C$ depending only on $M$), then, we easily get
that $c\|U-\hat U\|_{H(G)}\leq C\|G-\hat G\|_{C^1(\R)}$, where we have used that we have a uniform bound
of $\|\hat U\|_{H(G)}$ for all $G\in A(M)$, which is easily obtained from the variational formulation (\ref{weak-form-hat G}) and the fact that $\hat G\in A(M)$. This shows (\ref{conv-h1}) and we conclude the proof of the result. \end{proof}


\begin{thebibliography}{99}
\bibitem{ABMG} Y. Amirat, O. Bodart, U. de Maio, A. Gaudiello, ``Asymptotic Approximation 
of the solution of the Laplace equation in a domain with highly oscillating boundary'', \emph{SIAM J. Math. Anal.} 
35, 1598-1616 (2004)

\bibitem{JA} J. M. Arrieta, \emph{Spectral properties of Schr\"odinger operators under 
perturbations of the domain}, Ph.D. Thesis, Georgia Institute of Technology, (1991)
\bibitem{AP} J. M. Arrieta and M. C. Pereira,  ``Elliptic problems in thin domains with highly oscillating
boundaries'', \emph{Bol. Soc. Esp. Mat. Apl.} 51, pp:17-24 (2010).
\bibitem{ACPS} J. M. Arrieta, A. N. Carvalho, M. C. Pereira and R. P. da Silva;
 ``Nonlinear parabolic problems in thin domains with a highly oscillatory boundary'', \emph{Submitted}.
\bibitem{BLP} A. Bensoussan, J. L. Lions and G. Papanicolaou, \emph{Asymptotic Analysis for Periodic Structures},
North-Holland Publishing Company (1978).

\bibitem{BCh} R. Brizzi, J.P. Chalot, ``Boundary homogenization and Neumann boundary problem'' 
\emph{Ricerce di Matematica} XLVI, 2 (1997) 341-387

\bibitem{BL} V. Burenkov, P.D. Lamberti, ``Spectral Stability of general non-negarive self-adjoint operators 
with applications to Neumann-type operators'', \emph{J. Differential Equations} 233 (2007), 345-379


\bibitem{CP} D. Cioranescu and J. Saint Jean Paulin; \emph{Homogenization of Reticulated Structures},
Springer Verlag (1980).

\bibitem{DP}  A. Damlamian, K. Pettersson, ``Homogenization of oscillating boundaries''	, 
 \emph{Discrete and Continuous Dynamical Systems} 23, (2009), 197-219

\bibitem{HR} J. K. Hale and G. Raugel, ``Reaction-diffusion equation on thin domains'',
\emph{J. Math. Pures and Appl.} (9) 71, no. 1, 33-95 (1992).
\bibitem{Raugel} G. Raugel;
\newblock {\em Dynamics of partial differential equations on thin domains}
\newblock in Dynamical systems (Montecatini Terme, 1994), 208-315, 
Lecture Notes in Math., 1609, Springer, Berlin, 1995. 

\bibitem{SP} E. S\'anchez-Palencia, {\em Non-Homogeneous Media and Vibration Theory}, Lecture Notes in Physics 127, Springer Verlag (1980)

\bibitem{T1} L. Tartar; \emph{Probl\`emmes d'homog\'en\'eisation dans les \'equations aux d\'eriv\'ees partielles}, Cours Peccot, Coll\`ege de France (1977).
\bibitem{T2} L. Tartar, ``Quelques remarques sur l'homeg\'en\'eisation'', {\em Function Analysis and Numerical
Analysis}, Proc. Japan-France Seminar 1976, ed. H. Fujita, Japanese Society for the Promotion of Science, 
468-482 (1978).

\end{thebibliography}
\end{document}